\documentclass[12pt,a4paper]{amsart}
\usepackage[usenames]{color}
\usepackage{times}
\usepackage{amsfonts}
\usepackage{amsthm}
\usepackage{amsmath}
\usepackage{psfrag}

\usepackage{times}
\usepackage{latexsym,color}
\usepackage{amssymb}
\usepackage{graphicx}

\usepackage{epsfig}

\usepackage{tikz}
\usepackage[all]{xy}
 \usepackage[hang]{caption}
 \usetikzlibrary{snakes}
\usetikzlibrary{arrows}

\advance \topmargin by -\headheight
\advance \topmargin by -\headsep
\evensidemargin \oddsidemargin
\newcommand{\GQ}{{\mathbb Q}}

\newcommand{\GR}{{\mathbb R}}
\newcommand{\GZ}{{\mathbb Z}}

\newcommand{\GN}{{\mathbb N}}

\usepackage{lineno}
\linenumbers

\newcommand{\gN}{\mbox{\rm \scriptsize I$\!$N}}
\newcommand{\gZ}{\mbox{\rm \scriptsize Z$\!\!$Z}}

\newcommand{\ab}{|}

\newtheorem{theorem}{Theorem}
\newtheorem{lemma}[theorem]{Lemma}
\newtheorem{corollary}[theorem]{Corollary}
\newtheorem{proposition}[theorem]{Proposition}
\newtheorem{question}{Question}
\newtheorem{example}{Example}
\newtheorem{definition}{Definition}
\newtheorem{remark}{Remark}
\newtheorem{conj}{Conjecture}

\begin{document}\nolinenumbers

\title[]{Languages of general interval exchange transformations}

\author[S. Ferenczi]{S\'ebastien Ferenczi}
\address{Aix Marseille Universit\'e, CNRS, Centrale Marseille, Institut de Math\' ematiques de Marseille, I2M - UMR 7373\\13453 Marseille, France.}
\email{ssferenczi@gmail.com}

\author[P. Hubert]{Pascal Hubert}
\address{Aix Marseille Universit\'e, CNRS, Centrale Marseille, Institut de Math\' ematiques de Marseille, I2M - UMR 7373\\13453 Marseille, France}
\email{hubert.pascal@gmail.com}

\author[L.Q. Zamboni]{Luca Q. Zamboni}
\address{Institut Camille Jordan\\
Universit\'e Claude Bernard Lyon 1\\
43 boulevard du 11 novembre 1918\\
F69622 Villeurbanne Cedex
(France)}
\email{zamboni@math.univ-lyon1.fr}

\date{November 24, 2022}

\begin{abstract}
 The languages generated by interval exchange transformations have been characterized by Ferenczi-Zamboni (2008) and Belov-Cernyatev (2010) under some extra conditions on the system. Lifting these conditions leads us to consider successively natural codings of standard interval exchange transformations, natural codings of affine interval exchange transformations, grouped codings of affine interval exchange transformations, and natural codings of generalized interval exchange transformations. We show that these four classes of languages are strictly increasing, and give necessary and/or sufficient (but not all equally explicit) combinatorial criteria to describe each of them. These work also, mutatis mutandis, for interval exchanges with flips.
 \end{abstract}
\maketitle

Interval exchange transformations were originally
introduced by Oseledec
\cite{ose}, following an idea of Arnold \cite{arn}, see also \cite{ks}; the unit interval is partitioned
into  subintervals which are  rearranged by piecewise translations.
These {\em  standard} exchange transformations have been successively generalized to:
\begin{itemize}
\item interval exchange transformations {\em with flips} \cite{nog1}, where some or all translations are replaced by affine maps of slope $-1$;
\item {\em affine} interval exchange transformations, where the translations are replaced by affine maps of arbitrary nonzero slope (if some are negative we have a flipped affine interval exchange transformation), thus the Lebesgue measure is not preserved anymore; these appear in the precursor \cite{arth}, the oldest published reference may be \cite{lev};
\item {\em generalized} interval exchange transformations, where the translations are replaced by any continuous monotone nonconstant maps (if some are decreasing we have a flipped generalized interval exchange transformation), see \cite{arth} \cite{mmy2};
\item {\em systems of piecewise isometries}, see \cite{glp}, where the intervals are moved by isometries, but they or their images may be non-disjoint, thus the  transformation may be non-injective; these include the {\em interval translation mappings} of \cite{bos}.
\end{itemize}

These generalizations have seen  a recent surge in activity (see \cite{mmy} \cite{mmy2} \cite{gha} \cite{ulgh}) primarily centered on the conjugacy problem between these different classes of maps; in this context,  standard and generalized interval exchange transformations are the extreme cases while affine interval exchange transformations constitute a fundamental middle step.

 It was Rauzy \cite{rau} who first suggested using
interval exchange transformations as a possible framework for
generalizing the well-known interaction between dynamical systems, arithmetics and word combinatorics of   \cite{hm} and  \cite{ch}. In these papers, the languages  which arise as {\em natural codings} of minimal rotations by two intervals, aka {\em Sturmian} languages, are identified with the  {\em uniformly recurrent} languages of {\em complexity} $n+1$.  Rauzy  asked which  languages arise as  natural codings
of interval exchange transformations  through their defining intervals (see Definition \ref{sy}), the Sturmian case corresponding to an exchange of two intervals.  Several partial answers to Rauzy's question were given for $3$-interval exchange transformations in \cite{did}  \cite{san} \cite{fhz2}. Then two fairly general answers appeared independently  for all numbers of intervals. In \cite{fz3}, see Theorem \ref{tido} below, two of the authors of the present paper solved the problem for the most studied case among standard unflipped interval exchange transformations, those satisfying a version of  Keane's {\em i.d.o.c. condition}, a condition defined in \cite{kea} where it is shown to be stronger (indeed, strictly stronger) than minimality.  This description makes use of two main ingredients: an {\em order condition} (see Definition \ref{dord}) which dictates the way {\em bispecial} words can be {\em resolved} (see Section \ref{ud}),  and the absence of {\em connections} (see Definition \ref{dcon}). In contrast, the criterion in  \cite{bc} uses a description of the evolution of {\em Rauzy graphs} which is somewhat cumbersome to state, but applies to all standard interval exchange transformations, including those with flips. \\

In the present paper, we aim to study the languages of interval exchange transformations in the most general framework. We prove first that for every transformation which is called an interval exchange, the language satisfies an  order condition (Proposition \ref{pord}), with a new definition (Definition \ref{dord}) in the flipped cases, allowing the  orders to change according to the number of flipped letters in each bispecial word. The order condition is intrinsically linked to the one-dimensional nature of interval exchange transformations, and constitutes a strong constraint on the possible ways of resolving the bispecials, thus those languages arising from interval exchange transformation form a small family, which becomes larger when all flips are allowed. The order condition fails in general for piecewise isometries (see Remark \ref{pis}), and the characterization of the associated languages is completely unknown.

When the order condition is satisfied, to which extent do we need additional conditions? The following three counter-examples will be described in Section \ref{sexa} below: they all satisfy an unflipped order condition.
\begin{itemize}
\item The  language of Example  \ref{fs1} is not a natural coding of  a standard interval exchange transformation, but is the  natural coding of an affine interval exchange transformation.
\item  The language of Example \ref{fs2} is not a natural coding of an affine interval exchange transformation, but is the  coding of an affine interval exchange transformation, where the slope may take more than one value inside the coding intervals, we call it  a {\em grouped coding} of an affine interval exchange transformation.
\item The  language of Example  \ref{mon} is not a natural or grouped coding of an affine interval exchange transformation, but is the  natural coding of a generalized interval exchange transformation.
\end{itemize}

Thus all these four families of languages are different, and we can characterize some cases by simple word-combinatorial conditions. Under an order condition
\begin{itemize}
\item by Theorem \ref{tg2}, {\em all}
 languages correspond to generalized interval exchange transformations,
 \item by Theorem \ref{tg1}, {\em recurrent}
languages correspond to standard interval exchange transformations,
\item additionnally, by Theorem \ref{tmin}, {\em  aperiodic uniformly recurrent languages}  correspond to minimal standard interval exchange transformations.
\end{itemize}

But for the affine case, the characterization of the languages is more complicated, and is very much linked to the extensive work of \cite{lev} \cite{lio}  \cite{cage} \cite{bhm} \cite{co} \cite{mmy}, who build affine interval exchange transformations with wandering intervals semi-conjugated to a given standard interval exchange transformation. Thus our criteria to characterize their codings use Birkhoff sums as in these papers: what we get, under an order condition, are necessary or sufficient criteria for a language to  correspond to an affine interval exchange transformation by a natural coding (Theorem \ref{paf} and Corollary \ref{caf}), and a not very explicit necessary and sufficient criterion for a language to  correspond to an affine interval exchange transformation by a grouped coding (Theorem \ref{tg3}). 

Then Section \ref{sexa} is devoted to the  examples distinguishing our families of languges. A still open problem is to find an aperiodic grouped coding of an affine interval exchange transformation which is not a natural coding of an affine interval exchange transformation; it is quite possible that these do not exist, and in that case the problems of characterizing the natural and  grouped codings of  affine interval exchange transformations might get a natural solution (Conjecture \ref{con1} and Question  \ref{con2}). Another problem is to build  a language which is a natural coding of a generalized interval exchange transformation, but which is not a grouped coding of an affine interval exchange transformation,  which reveals different behaviours for affine and generalized interval exchange transformations. In Example \ref{mon}, the first property comes from the fact that, by a blow-up process due to Denjoy \cite{den}, it is possible to build a generalized interval exchange transformation $T$, with a wandering interval, semi-conjugate to a given standard interval exchange transformation $T'$ (this is indeed the basis of the proof of our Theorem \ref{tg2}), but if $T'$ is conjugate to a rotation,  the Denjoy-Koksma inequality \cite{her} prevents such $T$ from being affine, even after subdividing the continuity intervals; we also prove that $T$ fails to satisfy some regularity condition called class $P$ in \cite{her}. Another way to get such examples, allowing $T'$ to be any non purely periodic interval exchange transformation (not only a rotation) is given in Theorem \ref{ttrok} using Rokhlin towers. By using the results of \cite{mayo}, we have also been able to build exmmples where $T'$ is an interval exchange transformation naturally defined from any one of two famous translation surfaces, the {\it Eierlegende Wollmilch Sau} and the {\it Ornithorynque} (Section \ref{pascal}) but we could prove only that they are not natural codings of an affine interval exchange transformation. 

All our results extend without further difficulty to flipped interval exchange transformation by using flipped order conditions  (our counter-examples, however, are only written for the particular case where there are no flips).

\section{Languages}
\subsection{Usual definitions}\label{ud}

 Let $\mathcal A$ be a finite set called the {\em  alphabet}, its elements being {\em letters}.
 A {\em word} $w$ of {\em length} $n=|w|$ is $a_1a_2\cdots a_n$, with $a_i \in {\mathcal A}$. The {\em concatenation} of two words $w$ and $w'$ is denoted by $ww'$.\\
 
By a   language $L$ over $\mathcal A$ we mean  a {\it factorial extendable language}:
a collection of sets $(L_n)_{n\geq 0}$ where the only element of $L_0$ is the {\em empty word}, and where each $L_n$ for $n\geq 1$ consists of
words of length $n$,  such that
for each $v\in L_n$ there exists $a,b\in \mathcal A$ with $av,vb\in L_{n+1}$,
and
each $v\in L_{n+1}$ can be written in the form $v=au=u'b$ with $a,b\in \mathcal A$ and $u,u'\in
L_n.$  

 The {\em  complexity function} $p\,:\,\GN \rightarrow \GN$ is
defined by $p(n)=\# L_n$.

 A word $v=v_1...v_r$ is a {\em factor} of a word $w=w_1...w_s$ or an infinite sequence  $w=w_1w_2...$ if $v_1=w_i$, ...$v_r=w_{i+r-1}$.\\

 Let $W$ be a family of words or (one- or two-sided) infinite sequences. Whenever the set $L$ of  all the factors of the words or sequences in $W$ is a language (namely, is  factorial and extendable), we say that $L$ is {\em the language generated by  $W$} and denote it by $L(W)$.\\

 A language $L$ is
{\em recurrent} if for each $v\in L$ there exists  a nonempty $w\in L$, such that $vw$ ends with $v$.

A language $L$ is
{\em uniformly recurrent} if for each $v\in L$ there exists  $w\in L_n$ such that $v$
is a factor of each word $w\in L_n.$

A language $L$  is {\em aperiodic} if for  all nonempty words $w$ in $L$, there exists $n$ such that  $w^n$ is not in $L$.\\

  The {\em Rauzy graph} $G_n$ of a language $L$ is a directed graph whose vertex set consists of all words of length $n$ of $L$, with an edge from $w$ to $w'$ whenever $w=av$, $w'=vb$ for letters $a$ and $b$, and the word $avb$ is in $L$.\\

 For a word $w$ in $L$, we call {\em arrival set of $w$} and denote by $A(w)$ the set of all letters $x$ such that $xw$ is in $L$, and call {\em departure set of $w$} and denote by $D(w)$ the set of all letters $x$ such that $wx$ is in $L$.

 A word $w$ in $L$ is called {\em  right special}, resp. {\em
left
special} if $\# D(w)>1$, resp. $\# A(w)>1$. If $w\in L$ is both right special and
left special, then $w$ is called {\em  bispecial}. If $\# L_1>1$, the empty word $\varepsilon$
 is bispecial, with $A(\varepsilon)=D(\varepsilon)=L_1$.

 A bispecial word $w$ in $L$ is a {\em weak bispecial} if
 $\#\{awb \in L, a\in A(w), b\in D(w)\} < \#A(w)+\#D(w)-1$. 
 
 A bispecial word $w$ in $L$ is a {\em neutral bispecial} if
 $\#\{awb \in L, a\in A(w), b\in D(w)\} = \#A(w)+\#D(w)-1$. 
 
 A bispecial word $w$ in $L$ is a {\em strong  bispecial} if
 $\#\{awb \in L, a\in A(w), b\in D(w)\} > \#A(w)+\#D(w)-1$.

 To {\em resolve} a bispecial word $w$ is to find all words in $L$ of the form $awb$ for letters $a$ and $b$.
\\

The {\em symbolic dynamical system} associated to a  language $L$ is the two-sided shift $S$ acting on the subset $X_L$ of ${\mathcal A}^{\gZ}$ consisting of all bi-infinite sequences such that  $x_r\cdots x_{r+s-1}\in L_s$ for each $r$ and $s$, defined by $(Sx)_n=x_{n+1}$ for all $n\in\GZ$.

For a word $w=w_0\cdots w_t\in L$, the {\em cylinder} $[w]$ is the set $\{x\in X_L: x_0=w_0, \ldots , x_t=w_t\}$. \\

 In many papers, $(X_L, S)$ is the one-sided shift the subset $X_L$ of ${\mathcal A}^{\gN}$ consisting of all infinite sequences such that $x_r\cdots x_{r+s-1}$ is in $L_s$  for every $r,s$. In the present paper, we shall use two-sided sequences $x\in X_L$, but also their {\em  infinite suffixes} $(x_n, n\geq k)$ and {\em infinite prefixes} $(x_n, n\leq k)$.

\subsection{Languages  with order conditions}\label{lod}
\begin{definition}\label{dord}
A language $L$ on an alphabet $\mathcal A$ satisfies a {\em local order condition} if, for each bispecial word $w$, there exist two orders on $\mathcal A$, denoted  by $\leq_{A,w}$ and $\leq_{D,w}$, respectively $<_{A,w}$ and $<_{D,w}$ for the strict orders, such that
whenever  $awc$ and $bwd$ are in $L$ with letters $a\neq b$ and $c\neq d,$ then $a<_{A,w} b$ if and only if $c<_{D,w} d$.

A language $L$ on an alphabet $\mathcal A$ satisfies an $\mathcal F${\em -flipped order condition} for  a (possibly empty) subset $\mathcal F$ of $\mathcal A$ if there exist two orders on $\mathcal A$, denoted  by $\leq_A$ and $\leq_D$, respectively $<_A$ and $<_D$ for the strict orders, such that $L$ satisfies a local order condition where 
\begin{itemize}
\item the order $<_{A,w}$ is the same as $<_A$ for all  $w$,
\item $<_{D,w}$ is the same as  $<_D$ when the number of occurrences in $w$ of letters belonging to $\mathcal F$ is even, 
\item  $<_{D,w}$ is the reverse order of $<_D$ when the number of occurrences in $w$ of letters belonging to $\mathcal F$ is odd.
\end{itemize} \end{definition}

\begin{definition}\label{dcon}
In a language $L$, a {\em  locally strong bispecial word} is  a  bispecial word $w$ such that  there  exist nonempty subsets $A' \subset A(w)$, $D' \subset D(w)$ such that
 $\#\{awb \in L, a\in A', b\in D'\} > \#A'+\#D'-1$.

If a language $L$ on an alphabet $\mathcal A$ satisfies a local order condition, a bispecial word $w$ has a {\em connection} if there are letters $a<_{A,w}a'$, consecutive in the order $<_{A,w}$, letters $b<_{D,w}b'$, consecutive in the order $<_{D,w}$, such that $awb$ and  $a'wb'$ are in $L$, and neither $awb'$ nor $a'wb$ is in $L$.
\end{definition}

In this section, we state general combinatorial properties of languages satisfying some order conditions, or weaker properties, to be used in the next  sections. Note that some of these use the same ideas as Section 2.2 of \cite{fh4}, though the present context is both  more combinatorial and more general. Further properties with a more dynamical flavor will be studied at the beginning of Sections \ref{gai} and \ref{sai}.

\begin{lemma}\label{lsbs} A language $L$ which satisfies a local order condition  contains no   locally strong bispecial word, and thus no strong bispecial word. \end{lemma}
{\bf Proof}\\
Let $w$, $A'$, $D'$ be as in  Definition \ref{dcon}. If $\#A'$ or $\#D'$ is $1$, then the result is immediate. Assume $\#A'=\#D'=2$ and put $A'=\{a_1,a_2\}$, $D'=\{b_1,b_2\}$, $a_1<_A a_2$, $b_1<_D b_2.$ If $a_1wb_2$ and $a_2wb_1$ both belong to $L,$ this contradicts the local order condition. Now suppose we have proved the result for all $A'$ and $ D'$ such that $\#A'\leq p$ and $\#D'\leq q$; take now  $A'=\{a_1 <_{A,w} \cdots <_{A,w} a_p<_{A,w} a_{p+1}\}$, $D'=\{b_1 <_{D,w} \cdots <_{D,w} b_q\}$; then, let $r$ be the  largest integer between $1$ and $q$ such  that $a_pwb_r\in L.$  By the local order condition,  $a_iwb_j\in L$  for some $1\leq i\leq p$ only if $j\leq r$, thus by the induction hypothesis there are at most $p+r-1$ possible $a_iwb_j$, $1\leq i\leq p$, $1\leq j\leq q$. Again by the local order condition,  if $a_{p+1}wb_j\in L$ $1\leq j\leq q$,  then  $r\leq j\leq q$. Thus the number of possible $awb$ belonging to $L$ with $a\in A'$ and $b\in D'$ is at most $p+r-1+q+1$ and we have extended the induction hypothesis to $p+1$ and $q.$ A similar reasoning extends it to $p$ and $q+1$, thus the result holds for all possible $A'$ and $D'$. \qed\\

 For languages where each word has at most two right (resp. left) extensions, the absence of strong bispecial words, the absence of locally strong bispecial words, and a local order condition are all equivalent. In the general case, it is easy to find bispecials which are locally strong but not strong (suppose for example that the possible $xwy$ are $awa$, $awb$, $bwa$, $bwb$, $cwc$), and  a local order condition is stricter than the absence of locally strong bispecials.

\begin{example} Suppose $L$ is a language whose words of length $2$ are $ac$, $ad$, $ba$, $bc$,
$cb$, $cc$, $da$. Then the empty word is not a locally strong bispecial, yet $L$ does not satisfy any local order condition.

Indeed, assume to the contrary that $\varepsilon$ satisfies the order condition with respect to $<_{A, \varepsilon}$ and $<_{D, \varepsilon}$. By possibly reversing both orders, we can assume $c<_{A, \varepsilon} a.$
Since $\{cb,ac,cc,ad\}\subset L$  this implies
$b <_{D, \varepsilon} c<_{D, \varepsilon} d.$ This in turn implies $c <_{A, \varepsilon} b<_{A, \varepsilon} a$ (since $\{cb, bc, ad\}\subset L)$ which in turn implies $c <_{D, \varepsilon} a<_{D, \varepsilon} d$ (since $\{ba, cc, ad\}\subset L).$ Also $b <_{A, \varepsilon} d<_{A, \varepsilon} a$ because $\{bc, da, ad\}\subset L.$ Finally, as $ac$ and $da$ are in $L$, we get $a<_{D, \varepsilon} c$ which is a contradiction.

 We can then choose $L_3$ to be made with $acc$, $ada$, $bac$, $bad$, $bcb$, $bcc$, $cba$, $cbc$, $ccb$, $dac$, where each word has at most two left (resp. right) extensions and continue by resolving the bispecials so that they are all neutral. We get  a language without locally strong bispecials but not satisfying any local order condition. \end{example}
 
 \begin{lemma}\label{wbs} A language $L$ which has no strong bispecial factor has a finite number  of weak  bispecial factors. \end{lemma}
{\bf Proof}\\
Let $s(n)=p(n+1)-p(n)$. By Theorem 4.5.4 of \cite{cant}, see also \cite{casbel}, if $L$ has only weak or neutral bispecials $s(n+1)\leq s(n)$ for all $n$, and, if additionally $L$ has a weak bispecial of length $n$, $s(n+1)<s(n)$. This cannot happen infinitely many times as $p(n)$ is a growing function. \qed\\

Another notion aiming at generalizing the properties of some natural codings of interval exchange tranformations is a  {\it dendric} language, which is  defined in  \cite{bal} (under the name of tree sets): it turns out that a language is dendric iff it has neither locally strong bispecial words nor weak bispecial words, thus there is no inclusion relation between this family and the various order conditions. However, by Lemmas \ref{lsbs} and  \ref{wbs}, a language satisfying a local order condition is {\it ultimately dendric}.

\begin{lemma}\label{lcon} If a language $L$ satisfies a local order condition, a word $w$ is a weak bispecial iff it has a connection.  \end{lemma}
{\bf Proof}\\
Let $A(w)$, ordered by $<_{A,w}$, be $\{a_1,\ldots , a_p\}$, and $D(w)$, ordered by $<_{D,w}$, be $\{b_1,\ldots , b_q\}$. If there is no connection, then each  $D(a_kw)\cap D(a_{k+1}w)$ has at least one element, and indeed exactly one element by Lemma \ref{lsbs}: thus the number of possible $awb$, which is $\sum_{k=1}^p\#D(a_kw)$, is $q+p-1$, and $w$ is a neutral bispecial.

If there is a connection, then  for some $l$  $D(a_lw)\cap D(a_{l+1}w)$ is empty. By the local order condition, if $b_j$ is the maximal element (for $<_{D,w}$( of $D(a_lw)$, or equivalently of $D(a_1w)\cup \cdots D(a_lw)$, then
$b_{j+1}$ is the minimal element (for $<_{D,w}$( of $D(a_{l+1}w)\cup \cdots \cup  D(a_pw)$. By Lemma \ref{lsbs} applied first to $\{a_1,\ldots , a_l\}$ and $\{b_1,\ldots , b_j\}$, then to $\{a_{l+1},\ldots , a_p\}$ and $\{b_{j+1},\ldots , b_q\}$, we get that the number of possible $awb$ is at most $j+l-1+p-l+q-j-1=p+q-2$, thus $w$ is a weak bispecial.\qed\\

\begin{corollary}\label{compl} A language satisfying a local order condition has  complexity  $p(n)=kn+l$ for all $n$ large enough and with $0\leq k\leq \#\mathcal A-1.$ Moreover,  $k=\#\mathcal A-1$ if and only if $L$ has no connection, and in that case $l=1$. \end{corollary}
{\bf Proof}\\
This comes from Lemmas \ref{lsbs}, \ref{wbs}, \ref{lcon}, and
Theorem 4.5.4 of \cite{cant}.\qed\\

The following result is well known and could be deduced from Section 3.3 of \cite{casbel}, or Theorem 4.5.4 of \cite{cant}.

\begin{lemma}\label{orbls} If a language $L$ has no strong bispecial word,  the left (resp. right) special words are the prefixes (resp. suffixes) af a finite number of  infinite suffixes (resp. prefixes) of sequences of $X_L$.\end{lemma}
{\bf Proof}\\
Suppose $v$ is left special, and $va$ and $vb$ are left special for letters $a \neq b$, then both $A(va)$ and $A(vb)$ are strictly included in $A(v)$: indeed, otherwise, the possible $cvd$ in $L$ contain at least, for example $cva$ for all $c$ in $A(v)$, $c'vb$ for at least two different $c'$ in $A(v)$, and some $c_ivd_i$ for each $d_i$ in $D(v)\setminus \{a,b\}$, and thus $v$ is a strong bispecial word.

Starting from each left special letter $e$, we extend it to the right by one letter at a time; each time two different extensions of the same word are left special, the arrival sets decrease strictly; thus, after a finite number of bifurcations, each new left special extension has the previous one as a prefix, and all these are prefixes of a finite number of sequences $ex_1...x_n...$, which are infinite suffixes of sequences of $X_L$.

And a similar reasoning works  for the right specials. \qed\\

\begin{lemma}\label{deco} Let $L$ be a language satisfying a local order condition. If $G_n$ is  connected (as a nonoriented graph)  but  $G_{n+1}$ is not connected, then some   bispecial word in $L_n$ has a connection.\end{lemma}
{\bf Proof}\\
As $G_n$ is connected, if all pairs of vertices of $G_{n+1}$ which correspond to  two consecutive edges of $G_n$ belong to the same connected component of $G_{n+1}$, $G_{n+1}$ is connected.  Thus assume that for some $w\in L_n$ and letters $a$ and $b$, $aw$ and $wb$ are not  in the same connected component of $G_{n+1}$; then, in particular $awb$ is not in $L$. But as every word is extendable to the left and right, there exist $a'$ and $b'$ such that $a'wb$ and $awb'$ are in $L$, thus $a\neq a'$, $b\neq b'$, and $w$ is bispecial.

Suppose $w$ has no connection. If $a_1w$ is in some connected component $U$ of $G_{n+1}$, so is $wx$ for every $x$ in $D(a_1w)$, hence also $wx$  is in $U$ for  one element of  $D(a_2w)$ and so on, hence all the $a_iw$ and all $wx$ for $x$ in $D(w)$ are in $U$, and this contradicts the fact that $G_{n+1}$ is not connected.
\qed\\

Note that the terminology inherited from two different fields is somewhat counter-intuitive; a {\em connection} (so named for geometrical reasons) is necessary to  {\em disconnect} the Rauzy graphs. But it is not a sufficient condition, see examples in \cite{fh4}.\\

An order condition on $L$ allows us to define an order on the set $X_L$, which will be used in Section \ref{gai} below. It is well-defined because of the order condition.

\begin{definition}\label{otra} Suppose $L$ satisfies an $\mathcal F$-flipped order condition. Let $x$ and $y$ be in $X_L$. We define an order relation by $x<y$ if one or several of the following assertions is satisfied:
\begin{itemize}
\item $x_0<_D y_0$,
\item $x_i=y_i$ for $0\leq i \leq k$, $k>0$, $x_{k+1}<_Dy_{k+1}$,  an even number of the
$x_i, 0\leq i \leq k$, is in $\mathcal F$,
\item $x_i=y_i$ for $0\leq i \leq k$, $k>0$, $x_{k+1}>_Dy_{k+1}$,  an odd number of the
$x_i, 0\leq i \leq k$, is in $\mathcal F$,
\item $x_{-1}<_A y_{-1}$,
\item $x_i=y_i$ for $k\leq i \leq -1$, $k<-1$, $x_{k-1}<_Ay_{k-1}$,  an even number of the
$x_i, k\leq i \leq -1$, is in $\mathcal F$,
\item $x_i=y_i$ for $k\leq i \leq -1$, $k<-1$, $x_{k-1}>_Dy_{k-1}$,  an odd number of the
$x_i, k\leq i \leq -1$, is in $\mathcal F$.
\end{itemize}
\end{definition}

\begin{lemma}\label{lort} The order on $X_L$ is total, and every subset $Y$ of $X_L$ admits an upper bound  and a  lower bound in $X_L$. \end{lemma}
{\bf Proof}\\
For each $n$ we can find a word $y^{(n)}_m, -n\leq m\leq n$, such that $y^{(n)}$ is the  largest, resp. smallest, element of the set $\{y_m, -n\leq m\leq n, y\in Y\}$ for the order defined on such words by Definition \ref{otra}. Then $y^{(n+1)}_m=y^{(n)}_m$  for $-n\leq m\leq n$  and by definition of $X_L$ there is a sequence $y$ in $X_L$ such that $y^m=y^(n)_m$  for $-n\leq m\leq n$, which will be the desired  upper, resp. lower bound.  \qed\\

\section{Interval exchange transformations}\label{sie}
A generalized interval exchange
transformation is basically a bijection of
the unit interval to itself, which is continuous and monotone on a finite number of subintervals which partition the unit interval. However, such maps can be built, by using semi-open intervals of continuity, only when $T$ is increasing on all these intervals. In the {\em flipped} cases, $T$ will have to stay undefined on a finite set of points.

\begin{definition}\label{diet} Let $\mathcal A$ be a finite alphabet, $\mathcal F$ a possibly empty subset of
$\mathcal A$. An {\em $\mathcal F$-flipped generalized interval exchange transformation}  is  a map $T$ defined on a disjoint union of {\em open} intervals $I_e$, $e \in  \mathcal A$, such that the union of their  closures is $[0,1]$, continuous and (strictly) increasing on each $I_e$, $e \in  \mathcal F^c$, continuous and (strictly) decreasing on each $I_e$, $e \in  \mathcal F$, and such that the $TI_e$, $e\in \mathcal A$, are disjoint open intervals and the union of their closures is $[0,1]$.

The $I_e$, indexed in $\mathcal A$, are called the {\em defining intervals} of $T$.

If the restriction of $T$ to each $I_e$ is an affine map, $T$ is an {\em $\mathcal F$-flipped  affine} interval exchange transformation.

If the restriction of $T$ to each $I_e$ is an affine map of slope $\pm 1$, $T$ is an {\em $\mathcal F$-flipped standard} interval exchange transformation.

The endpoints of the $I_e$, resp. $TI_e$, excluding $0$ and $1$, will be denoted by $\gamma_i$, resp. $\beta_j$, for $i$, resp. $j$, taking $\#\mathcal A-1$ values. 
\end{definition}

In this definition, the defining intervals $I_e$ are not necessarily the intervals of continuity of $T$, as it may happen that $I_e$ and $I_f$ are adjacent and sent to adjacent intervals with the same flips in the same order, thus $T$ may be extended to a continuous map on their union. Thus in the present paper an interval exchange transformation $T$ is always supposed to be given {\em together with its defining intervals} as keeping the same $T$ but changing the defining intervals would change the coding.

 In all cases, $T$ is undefined on a finite number of points, namely the $\gamma_i$, $0$ and $1$. In the unflipped cases,  $T$ can be extended to $[0,1)$, resp. $(0,1]$, by including  in  each $I_e$ its left (resp. right) endpoint. Also in this case, as did  Keane in  Section 5 of \cite{kea} for unflipped standard interval exchange transformations and Arnoux in Chapter 1 of \cite{arth} (Definition 1.4 and Lemma 1.5) for unflipped generalized interval exchange transformations,  by carefully doubling the endpoints and their orbits, it is possible to  extend $T$ to an homeomorphism on the union of the unit interval and a Cantor set.

Throughout this paper, all defining intervals will be nonempty open intervals. In the unflipped cases, all our results and proofs remain valid if we define $T$ with semi-open intervals, open on the right (resp. left), or define $T$  on Keane's set (mentioned above) and use closed intervals.

\begin{definition}\label{wan}
Let $I$ be the subset (with countable complement) of $(0,1)$ where all $T^n$, $n\in \GN$, are defined. Every statement concerning $T^n$ will be tacitly assumed to be valid if it is valid on $I$.

$T$ is {\em minimal} if every orbit is dense.

A {\em wandering interval} is an interval $J$ for which $T^nJ$ is disjoint from $J$ for all $n>0$.
\end{definition}

\subsection{Codings}
\begin{definition}\label{sy} For a generalized $\mathcal F$-flipped interval exchange transformation $T$, coded by $\mathcal A$, its {\em natural coding} is  the  language
 $L(T)$ generated by all the {\em trajectories}, namely the
sequences
$(x_n, n\in\GZ) \in \mathcal A^{\GZ}$ where $x_{n}=e$ if ${T}^nx$ falls into
$I_e$, $e\in \mathcal A$.
\end{definition}

Thus we can look at the symbolic system associated to $L(T)$. Note that the set
 $X_{L(T)}$ is the closure in $\mathcal A^{\gZ}$ of the set of trajectories, for the product topology  defined by  the discrete topology on $\mathcal A$.

 \begin{example} A {\em Sturmian language} is the natural coding of the unflipped standard  interval exchange transformation $T$
sending $(0,1-\alpha)$  to $(\alpha,1)$ and $(1-\alpha, 1)$ to $(0,\alpha)$ for $\alpha$ irrational; $T$ is conjugate to a rotation of angle $\alpha$ on the $1$-torus. \end{example}

\begin{remark}\label{impp} The elements of $X_{L(T)}$ which are not actual trajectories of points under $T$ are called {\em improper}  trajectories.  For each improper trajectory $u$ in $X_{L(T)}$, either there exists an endpoint $\gamma$ of a defining interval $I_e$, and an integer $l\leq 0$ such that $u$ is the limit  of the trajectories of $x$ when $x$ tends to $T^l\gamma$, either from the left or from the right, or there exists an endpoint $\beta$ of a $TI_e$, and an integer $l\geq 0$ such that $u$ is the limit  of the trajectories of $x$ when $x$ tends to $T^l\beta$, either from the left or from the right.

In the unflipped cases, if we define $T$ on Keane's set (as mentioned after Definition \ref{diet}), each element of $X_{L(T)}$ is an actual trajectory for $T$; if we use semi-open intervals,  each element of $X_{L(T)}$ is an actual trajectory, either for $T$ defined and coded with semi-open intervals $[a,b)$, or  for $T$ defined and coded with semi-open intervals $(a,b]$, in the same way as Sturmian trajectories in \cite{hm}.\end{remark}

We shall also consider slightly more general codings, by merging  into intervals $\tilde I_e$ some adjacent intervals $I_e$ whose images by $T$  are also adjacent, and flipped in the same way. This is equivalent to taking the natural coding of another interval exchange transformation $\tilde T$, but when  $T$ is affine, if we define $\tilde T$ by the intervals $\tilde I_e$ it will not necessarily be affine by  our definition, as the slope is not constant on its defining intervals, see Example \ref{fs2} below. Thus we define

\begin{definition}\label{gc} A language $L$ is a {\em grouped coding} of an $\mathcal F$-flipped affine interval exchange transformation $T$ if there exist intervals $\tilde I_e$, $e\in\tilde{\mathcal A}$ such that
\begin{itemize}
\item each $\tilde I_e$ is an open interval, and a disjoint  union (plus intermediate endpoints) of defining intervals of $T$,
\item $T$ can be extended by continuity to a continuous monotone map on each $\tilde I_e$,
\item $L$ is the coding of $T$ by the $\tilde I_e$, that is the language generated by the trajectories $(x_n, n\in \GZ) \in \tilde{\mathcal A}^{\GZ}$ where $x_{n}=e$ if ${T}^nx$ falls into
$\tilde I_e$, $e\in \tilde{\mathcal A}$.
\end{itemize}
\end{definition}

\subsection{Orders}
 \begin{definition}\label{diord} A generalized interval exchange transformation defines two orders on $\mathcal  A$:
\begin{itemize}
\item $e<_D f$ whenever the interval $I_e$ is strictly to the left of the interval $I_f$,
\item $e<_A f$ whenever the interval $TI_e$ is strictly to the left of the interval $TI_f$.
\end{itemize}
\end{definition}

These orders correspond to the two permutations  used by Kerckhoff \cite{ker} to define
 standard interval exchange transformations: the unit interval is partitioned into semi-open intervals
which are numbered from $1$ to $k$, ordered according to a permutation  $\pi_0$ and then rearranged according to another permutation  $\pi_1$; in more classical definitions, there is only one permutation $\pi$, which corresponds to $\pi_1$ while $\pi_0=Id$; note that sometimes the orderings are by $\pi_0^{-1}$ and $\pi_1^{-1}$.

\begin{proposition}\label{pord}
Let $T$ be a generalized $\mathcal F$-flipped interval exchange transformation, for $\mathcal F$ a subset of an alphabet $\mathcal A$. Then its language $L(T)$ satisfies an $\mathcal F$-flipped order condition.\end{proposition}

{\bf Proof}\\
 Let $<_A$ and $<_D$ be the orders from Definition \ref{diord}.

Let $w=w_0\cdots w_{r-1}$ be a bispecial word in $L(T).$. The cylinder
$[w]$ is the interval $I_{w_0}\cap T^{-1}I_{w_1}\cap \cdots \cap T^{-r+1}I_{w_{r-1}}$; after deleting a finite number of points, this interval is partitioned into the $T[aw]=TI_a \cap [w]$, $a\in A(w)$. As
$T[aw] \subset TI_a$ for all $a\in A(w)$, $T[aw]$ is (strictly) left of $T[a'w]$ whenever $a\leq_A a'$.

Similarly, after deleting a finite number of points, $[w]$ is also partitioned into the $[wb]=I_{w_0}\cap T^{-1}I_{w_1}\cap \cdots \cap T^{-r+1}I_{w_{r-1}}\cap T^{-r}I_b\subset T^{-1}I_{w_1} \cap \cdots \cap T^{-r+1}I_{w_{r-1}}\cap T^{-r}I_b$.

The order (from left to right) between two intervals $[b]=I_b$ and $[b']=I_{b'}$ is the order $<_D$. The order between $TI_{w_{r-1}} \cap I_b\subset I_b$ and
$TI_{w_{r-1}} \cap I_{b'}\subset I_{b'}$ is also the order $<_D$. The order between  $I_{w_{r-1}} \cap T^{-1}I_b$ and
$I_{w_{r-1}} \cap T^{-1}I_{b'}$ is the same order if $T$ is increasing on $I_{w_{r-1}}$ or equivalently $w_{r-1} \in\mathcal F^c$, the opposite order otherwise.
Similarly the order between $TI_{w_{r-2}} \cap I_{w_{r-1}} \cap T^{-1}I_b\subset I_{w_{r-1}} \cap T^{-1}I_b$ and
$TI_{w_{r-2}} \cap I_{w_{r-1}} \cap T^{-1}I_{b'} \subset I_{w_{r-1}} \cap T^{-1}I_{b'}$ is the same as the last one, thus the order between $I_{w_{r-2}} \cap T^{-1}I_{w_{r-1}} \cap T^{-2}I_b$ and
$I_{w_{r-2}} \cap T^{-1}I_{w_{r-1}} \cap T^{-2}I_{b'}$
is the same order as the last one  if $T$ is increasing on $I_{w_{r-2}}$ or equivalently $w_{r-2} \in
\mathcal F^c$, the opposite order otherwise. And so on, finally we get that the order between
$[wb]$ and $[wb']$ is exactly the order $<_{D,w}$ defined in the second item of Definition \ref{dord}.

Thus the intervals $T[aw]$ are ordered from left to right by the order $<_A$, while the $[wb]$ are ordered from left to right by the order $<_{D,w}$.
The word $awc$ exists whenever ${T}[aw]\cap [wc]$ is nonempty, similarly for $bwd$:  as in  Figure 1, we get that if $awc$ and $bwd$ exist, $a\neq c$, $b\neq d$, then $a<_Ac$ iff $b<_{D,w}d$. Note that this is true also for the empty word, which is bispecial and for which no letter  is in $\mathcal F$. Thus $L(T)$ does satisfy the $\mathcal F$-flipped order condition defined by $\leq_A$ and $\leq_D$. \qed\\

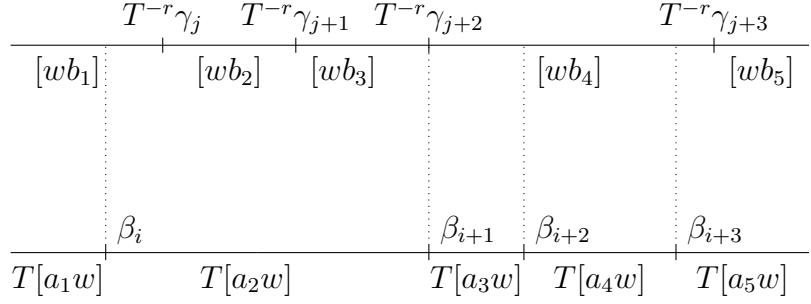
\begin{figure}
\begin{center}
\begin{tikzpicture}[scale = 5]

\draw (0,0.55)node[above]{};
\draw (.4,0.55)node[above]{$T^{-r}\gamma_j$};
\draw (.75,0.55)node[above]{$T^{-r}\gamma_{j+1}$};
  of an interval exchange transformation
\draw (1.1,0.55)node[above]{$T^{-r}\gamma_{j+2}$};

\draw (1.85,0.55)node[above]{$T^{-r}\gamma_{j+3}$};
\draw (2.1,0.55)node[above]{};

\draw(.4,.53)--(.4,.57);
\draw(.75,.53)--(.75,.57);

\draw(1.1,.53)--(1.1,.57);
\draw(1.85,.53)--(1.85,.57);

\draw(0,.55)--(.4,.55);
\draw(.4,.55)--(.75,.55);
\draw(.75,.55)--(1.1,.55);
\draw(1.1,.55)--(1.5,.55);
\draw(1.85,.55)--(2.1,.55);

\draw(1.5,.55)--(1.85,.55);

\draw (.15,0.55)node[below]{$[wb_1]$};
\draw (.575,0.55)node[below]{$[wb_2]$};
\draw (.875,0.55)node[below]{$[wb_3]$};
  \draw (1.47,.55)node[below]{$[wb_4]$};

\draw (1.975,0.55)node[below]{$[wb_5]$};

\draw (0,0)node[above]{};
\draw (.25,.06)node[above, right]{$\beta_{i}$};

\draw (1.1,.06)node[above, right]{$\beta_{i+1}$};
\draw (1.35,.06)node[above, right]{$\beta_{i+2}$};
\draw (1.75,.06)node[above, right]{$\beta_{i+3}$};
\draw (2.1,0)node[above]{};

\draw(0,0)--(.25,0);
\draw(.25,0)--(.65,0);
\draw(.65,0)--(1.1,0);
\draw(1.1,0)--(1.35,0);
\draw(1.35,0)--(1.75,0);
\draw(1.75,0)--(2.1,0);

\draw[dotted](.25,.02)--(.25,.55);

\draw[dotted](1.1,.02)--(1.1,.55);
\draw[dotted](1.35,.02)--(1.35,.55);
\draw[dotted](1.75,.02)--(1.75,.55);

\draw(.25,-.02)--(.25,.02);

\draw(1.1,-.02)--(1.1,.02);
\draw(1.35,-.02)--(1.35,.02);
\draw(1.75,-.02)--(1.75,.02);

\draw (.125,0)node[below]{$T[a_1w]$};
\draw (.62,0)node[below]{$T[a_2w]$};
  \draw (1.225,0)node[below]{$T[a_3w]$};
\draw (1.55,0)node[below]{$T[a_4w]$};
\draw (1.925,0)node[below]{$T[a_5w]$};

\end{tikzpicture}
\caption{A bispecial interval}
\end{center}
\end{figure}

\begin{remark}\label{ordint}By the same  reasoning as in  Proposition \ref{pord}, we check that if $(x_{n}, n\in \GZ)$ in $X_{L(T)}$ is the (actual or improper) trajectory of some point $x$, $(y_n, n\in \GZ)$ in $X_{L(T)}$ is the (actual or improper) trajectory of some point $y$, $(x_{n}, n\in \GZ)\leq (y_{n}, n\in \GZ)$ for the order of Definition \ref{otra} if and only if $x\leq y$
in the natural order on $(0,1)$.\end{remark}

\begin{remark}\label{pis} The language of an {\em interval translation mapping} \cite{bos} does not necessarily satisfy an $\mathcal F$-flipped order condition: it is possible that $TI_a$ intersects $TI_b$ for $a\neq b$, and, if $TI_a\cap TI_b$ intersects both $I_c$ and $I_d$, $c\neq d$, then the empty word is a locally strong bispecial, which contradicts the local order condition by Lemma \ref{lsbs}.\end{remark}

 \begin{question} How to characterize the   natural codings of systems of piecewise isometries,  resp.  interval translation mappings ?\end{question}
 
 Even  for the latter family, almost nothing is known. Their complexity is the object of a question of Boshernitzan; it was proved to be linear in some particular cases by Cassaigne (unpublished).

\section{Standard interval exchange transformations}
We continue the study   of languages satisfying  order conditions, or weaker properties, by looking at {\em recurrence}. Recurrence and uniform recurrence of a language are defined in Section \ref{ud} above. A different concept is the  recurrence of  infinite sequences, which is a property of the associated symbolic dynamical system. 

\begin{definition} A bi-infinite sequence $x$ in ${\mathcal A}^{\GZ}$, or an infinite suffix of it, is {\em right recurrent} if any factor $w$
of $x$ is equal to $x_{n_k}...x_{n_k+t}$ for  a sequence $n_k$ tending to $+\infty$.

A bi-infinite sequence $x$ in ${\mathcal A}^{\GZ}$, or an infinite prefix of it, is {\em left recurrent} if any factor $w$
of $x$ is equal to $x_{n_k}...x_{n_k+t}$ for  a sequence $n_k$ tending to $-\infty$.

A bi-infinite sequence $x$ in ${\mathcal A}^{\GZ}$ is {\em recurrent} if it is both left and right recurrent.
\end{definition}

\begin{lemma}\label{recu} Suppose $L$ has no strong bispecial word.

If an infinite suffix of a sequence in $X_L$ is right recurrent, it generates a uniformly recurrent language.

If an infinite prefix of a sequence in $X_L$ is left recurrent, it generates a uniformly recurrent language.

If a sequence in $X_L$ is left or  right recurrent, it is  recurrent and generates a uniformly recurrent language.\end{lemma}
{\bf Proof}\\
 Let  $w$ be any word in $L$; we look at the  possible {\em return words} of $w$, namely words $v\neq  w$ having $w$ as a prefix and as a suffix, but as no other subword. Extending  $w$ from the left, before seeing $w$ as a suffix,  there are several possible paths; bifurcations correspond to right special words of the form $wv$; by Lemma \ref{orbls} these correspond to suffixes of $K$ infinite prefixes of sequences in $X_L$. If there are more than $K$ bifurcations, some $wv_1$ and $wv_2$ are different suffixes of the same infinite prefix; thus there is an occurrence of $w$ in the shortest suffix which is not the initial occurrence, thus we have already returned to $w$ before the bifurcation; hence we cannot see more than $K$ bifurcations, with at most $K'$ choices  for each one,  before returning to $w$. Thus there are at most
$K'K$ return words of $w$.\\

If $x$ is right recurrent, as any prefix of $x$ occurs  further to the right, the set of its factors is a language $L'$. For each factor $w$ of $x$, we see it infinitely many times, with at most $KK'$ ways of going from one occurrence  to the next one. Thus we see $w$ in $x$ infinitely many times with bounded gaps, thus $w$ occurs in any
long enough factor of $x$, and $L'$ is uniformly recurrent. And the same works for left recurrence.

Now, suppose for example $x\in X_L$ is right recurrent. Then it is also left recurrent: otherwise there exists a factor $w$ of $x$ which does not occur in $(x_n, n\leq k)$; thus there exists arbitrarily long factors $v_{(m)}$ of $x$ with no occurrence of $w$; but then any infinite suffix of $x$ contains every $v_{(m)}$, and thus cannot generate an uniformly recurrent language.
Hence any factor $w$ occurs in $x$ at infinitely many places to the right and to the left, and with bounded gaps, thus the language generated by $x$ is uniformly recurrent.\qed\\

It will result from Proposition \ref{cmin} that, under a local order condition, the recurrence of $L$ implies the recurrence of every sequence in $X_L$; here we prove first a weaker statement.

\begin{lemma}\label{wapd} Suppose $L$ has no locally strong bispecial, and is recurrent. Then if $w'w^n$ is in $L$ for all $n \in \GN$, with $|w'|\leq |w|$, then $w'$ is a suffix of $w$, and if $w^nw'$ is in $L$ for all $n$, then $w'$ is a prefix of $w$. \end{lemma}
{\bf Proof}\\ Suppose first $w'w$ is  a non-initial factor of some $w'w^m$;  then either $w'$ occurs as a suffix of $w$, or we have $w=w_1...w_l=w_{d+1}...w_lw_1...w_d$ for some $0<d<l$; in the latter case, we get that $w=vv...v$ for a factor $v$ of length $p=d\wedge l$, and that $w'$ is a factor of $ww$ ending at index $d$ or $l-d$, and this factor is also a suffix of $w$. 

Suppose $w'w^n$ is in $L$ for all $n$, and $w'$ is not a suffix of $w$. By recurrence, for each $n\geq 1$ $w'w^n$ must occur after $w'w^n$. As it cannot occur as a non-initial factor of some $w'w^N$,  then for each $n\geq 1$ there is  a right special word which is a prefix of some $w'w^N$ and has  $w'w^n$ as a prefix; for infinitely many $p_n$, it will be of the form $w'w^{p_n}w_1$, for a fixed prefix $w_1$ of $w$, and followed by  two fixed letters $b$ and $b'$. If $w_2$ is the longest common suffix of $w$ and $w'$, we get four different words of the form $aw_1w^{p_n}w_2b$, and thus a locally strong bispecial. And similarly for the second assumption. \qed\\

\begin{lemma}\label{minm} A finite or countable union of uniformly recurrent languages satisfies the {\em measure condition}: there exists  an invariant probability measure $\mu$ on the symbolic system $(X_L,S)$ associated to $L$ such that $\mu [w]>0$ for any $w\in L$. \end{lemma}
{\bf Proof}\\
Assume $L$ is a uniformly recurrent language and let $\mu$ be any invariant probability measure on the symbolic system $(X_L,S)$ associated to $L.$  Such a measure $\mu$ can be constructed as in \cite{bos2}. Fix $w$ in $L.$ Then, for each positive integer $N$ there is at least one word $w_1$ of length $N$ with $\mu [w_1]>0$, and, as $w$ is a factor of $w_1$ if $N$ is large enough, $[w_1]$ is included in $\cup_{i=0}^NS^i[w]$ from which it follows that  $\mu [w]>0$.

If $L$ is a finite or countable union of uniformly recurrent languages $L_i$, then each $L_i$ satisfies the measure condition for some measure $\mu_i$, and hence $L$ will satisfy the measure condition for any average of the measures $\mu_i.$ \qed\\

\begin{proposition}\label{cmin} Let $L$ be a recurrent language satisfying a local order condition. Then $L$ is a finite union of uniformly recurrent languages. \end{proposition}
{\bf Proof}\\
 The Rauzy graph $G_1$ has at most $k$ connected components. The language $L$ restricted to any connected component of any of its Rauzy graphs $G_n$ still satisfies the local order condition. We apply repeatedly Lemma \ref{deco} to any of these languages: if $G_{n+1}$ has more connected components than  $G_n$, then there exists a connected component $C_n$ of $G_n$ such that the Rauzy graph  $C_{n+1}$ made with the edges of $C_n$ is not connected. Thus there is an edge missing in $C_{n+1}$ as in Lemma \ref{lcon}, and it cannot be elsewhere in $G_{n+1}$, thus this edge is missing in $G_{n+1}$ and there is  a weak bispecial word in $L$ by the reasoning of Lemma \ref{lcon}. By Lemma  \ref{wbs} this happens a finite number of times, thus for $n$ large enough $G_n$ has at most $N$ connected components, and $L$ is a union of $N$ languages for which all Rauzy graphs are connected. \\

Suppose there is a word $w$ in $L$ for which no left extension is left special. Then, before $w=w_0\cdots w_l$, we can see only one sequence of letters,  denoted by $w_{n}$, $n<0$.  As $L$ is recurrent,  $w$ must occur as some $w_s...w_t$, $s<t<0$, thus the sequence $(w_n, n<0)$ begins with a periodic infinite prefix $\ldots vvv$. For any $m$ large enough, this gives a loop in the Rauzy graph $G_m$; it is impossible both to enter this loop from any other point in $G_n$ and  to exit from it to any other point in $G_n$, otherwise this contradicts Lemma \ref{wapd}. Thus the language generated by $\ldots vvv$ corresponds to a full connected component of $G_m$ for all $m$ large enough, thus there are at most $N$ such languages. Thus  $L$ is the union of at most $N$  languages  satisfying the local order condition, where each language either is generated by the $v^n$, $n\geq 0$, for a word $v$,  or contains only   words which can be extended to a left special word.

Let $L'$ be one of these last languages, if they exist.  Let $G'_n$ be its Rauzy graphs. Because there is no strong bispecial word in $L'$, by  Lemma \ref{orbls} all left special words in $L'$ are prefixes of  $K$ infinite suffixes $W_i$; as every word of $L'$ can be extended to a left special, all the words of $L'$ are the factors of  the $W_i$.

Let $w$ be a word of $L'$: $w$ must be a suffix of infinitely many words of $L'$, because every word of $L'$ can be extended infinitely many times to the left; thus $w$ occurs infinitely often in at least one $W_i$; hence all the prefixes of $W_i$ occur infinitely often in at least one $W_j$, and if $j\neq i$, then we can drop $W_i$ and use only the other $W_j$ to generate $L'$. Thus, after a renumbering of the $W_i$,  $L'$ is the union of the languages $L(W_i)$, $1\leq i\leq q\leq K$, where the $W_i$ are right recurrent infinite suffixes. \\

Then by Lemma \ref{recu} each $L(W_i)$ is uniformly recurrent, while any of  language generated by the $v^n$, $n\geq 0$, is also uniformly recurrent.\qed\\

We can now prove our first main theorem, which we guess might also be deduced from understanding and rewriting the criterion of \cite{bc}.

\begin{theorem}\label{tg1}
For a language $L$ on an alphabet $\mathcal A$ and
$\mathcal F\subset \mathcal A$, the following are equivalent:
\begin{itemize}
\item{$(i)$} $L$ satisfies an $\mathcal F$-flipped order condition and is recurrent;
\item{$(ii)$} $L$ is the language of a standard $\mathcal F$-flipped interval exchange transformation;
\item{$(iii)$} $L$ is the language of a generalized
 $\mathcal F$-flipped  interval exchange transformation without wandering intervals. \end{itemize}
\end{theorem}
{\bf Proof}\\
If $L$ satisfies $(i)$, by Proposition \ref{cmin} and
Lemma \ref{minm}, $L$ satisfies the measure condition of Lemma \ref{minm}.
We prove now that it implies $(ii)$; this is indeed the same reasoning as in Proposition 4 of \cite{fz3}, modified to take into account flips and connections.\\

Let $\mu$ be as in the measure condition. Let  $\mathcal A=L_1$, and let $T$ be the standard $\mathcal F$-flipped  interval exchange transformation defined by intervals $I_e$, $e\in \mathcal A$, of respective length $\mu[e]$ ordered from left to right by the order $<_D$, while the $TI_e$, $e\in \mathcal A$, of respective length $\mu[e]$,  are ordered from left to right by the order $<_A$, and $T$ has slope $1$ on the $I_e$, $e\in \mathcal F^c$, and slope $-1$ on the $I_e$, $e\in \mathcal F$.

  Let $L'=L({T})$, $\mu'$ the shift-invariant measure on $X_{L'}$ defined by the Lebesgue measure on $(0,1)$.

We prove inductively on $n$ that $w \in L_n$ if and only if $w\in L'_n$,  and for these words $\mu [w]=\mu  '[w]$. This is true for $n=1$ by our choice of $T$. We denote by $A_L$, $A_{L'}$, $ D_L$, $D_{L'}$ the arrival and departure sets in $L$ and $L'$.

Let the hypothesis be proved for $n$, and take a word $w$ in $L$ of length $n-1$ (possibly the empty one if $n=1$(.  Because of the induction hypothesis, $A_L(w)=A_{L'}(w)$, $D_L(w)=D_{L'}(w)$, as they depend only on $L_n$ and $L'_n$.

If $w$ is not right special, $D_L(w)$ is made with a single letter $a$, and every word $xw$ is always followed by $a$, thus $xwa$ is in $L$ if and only if  $xw$ is in $L$, with $\mu[xwa]=\mu[xw]$; because of our induction hypothesis, all this remains true with ``prime" signs added. Thus our hypothesis is carried over to all $awb$, for $a$ and $b$ letters. And similarly if $w$ is not left special.\\

Suppose now that $w$ is bispecial in $L$ (hence also in $L'$).

Let $A_L(w)=\{a_1, \ldots ,a_p\}$, in the $<_A$ order, and $D_L(w)=\{b_1, \ldots ,b_q\}$  in the $<_{D,w}$ order defined in the second item of Definition \ref{dord}. Then what we shall show is that the situation is analogous to Figure 1 above.

Starting from the left, $D_L(a_1w)$ must contain $b_1$, otherwise $D_L(a_1w)$ contains some  $b>_D b_1$, while  $A_L(wb_1)$ contains some  $a>_A a_1$, $a_1wb$ and $awb_1$ exist, this contradicts the order condition. Then $a_1wb_1$ is in $L$, and
\begin{itemize} \item if $\mu[a_1w]\leq  \mu[wb_1]$, then, as $\mu$ gives positive measure to every cylinder,  $D_L(a_1w)$ is reduced to the element $b_1$ and thus $\mu[a_1wb_1]=\mu[a_1w]$,
\item if $\mu[a_1w]> \mu[wb_1]$, then $D_L(a_1w)$ contains also another element than $b_1$, thus it must  contains $b_2$, otherwise $a_1wb$ and $awb_2$ exist for some $b>_{D,w} b_2$ and $a>_A a_1$; hence $A(wb_2)$ contains $a_1$,  $A(wb_1)$ cannot contain $a>_A a_1$, thus is reduced to the element $a_1$  and $\mu[a_1wb_1]=\mu[wb_1]$.
\end{itemize}

Thus $\mu[a_1wb_1]=\min (\mu[a_1w], \mu[wb_1])>0$; similarly, because of the definition of $L'$, $a_1wb_1$ is in $L'$, and, using our induction hypothesis, we get
$$\mu'[a_1wb_1]=\min (\mu'[a_1w], \mu'[wb_1])=\min (\mu[a_1w], \mu[wb_1])=\mu[a_1wb_1].$$

Suppose we have proved that $a_rwb_s$ is in $L$ iff it is in $L'$, and $\mu[a_rwb_s]=\mu'[a_rwb_s]$, for every $r\leq i$, and $s\leq j$; suppose that $a_iwb_j$ is in $L$ and $L'$, with $i<p$ or $j<q$.

 Suppose first that  $\sum_{r=1}^i\mu[a_rw]> \sum_{s=1}^j\mu[wb_s]$; then $j<q$ and $\cup_{r=1}^iD_L(a_rw)$ contains strictly $\{b_1,\ldots , b_j\}$, thus by the order condition $D_L(a_iw)$ contains $b_{j+1}$, thus $a_iwb_{j+1}$ is in $L$ while $a_{i+1}wb_{j}$ is not in $L$. And we have again two cases: \begin{itemize}
\item  if $\sum_{r=1}^i\mu[a_rw]> \sum_{s=1}^{j+1}\mu[wb_s]$, then by the order condition $A_L(wb_{j+1})$ is reduced to the element $a_i$, thus $\mu[a_iwb_{j+1}]=\mu[wb_{j+1}]$,
\item  if $\sum_{r=1}^i\mu[a_rw]< \sum_{s=1}^{j+1}\mu[wb_s]$, then by the order condition $D_L(a_iw)$ contains (at least) the two elements $b_j$ and $b_{j+1}$, and either
$\mu[a_iwb_{j+1}]=\sum_{r=1}^i\mu[a_rw]- \sum_{s=1}^j\mu[wb_s]$, or $\mu[a_iwb_{j+1}]= \sum_{s=1}^{j+1}\mu[wb_s]- \sum_{r=1}^i \mu[a_rw]$, according to which of these two quantities is positive.
\end{itemize}
 But then  $\sum_{r=1}^i\mu'[a_rw]> \sum_{s=1}^j\mu'[wb_s]$, and the same analysis applies to $L'$, thus we get the same conclusions with ``prime" signs added; as all estimates depend only on properties of words of length $n$, we use our induction hypothesis to get $\mu'[a_iwb_{j+1}]=\mu[a_iwb_{j+1}]$, and our induction hypothesis is carried to $i$ and $j+1$.

 In the opposite case where
 $\sum_{r=1}^i\mu[a_rw]< \sum_{s=1}^j\mu[wb_s]$, then $i<p$, $a_iwb_{j+1}$ is not in $L$ while $a_{i+1}wb_{j}$ is in $L$, and a similar reasoning applies, carrying our induction hypothesis to  $i+1$ and $j$.

If  $\sum_{r=1}^i\mu[a_rw]= \sum_{s=1}^j\mu[wb_s]$, then there is a  connection. Then $i<p$, $j<q$,  by the reasoning above $\mu[a_iwb_{j+1}]=\mu[a_{i+1}wb_j]=\mu'[a_iwb_{j+1}]=\mu'[a_{i+1}wb_j]=0$ and the words
$a_iwb_{j+1}$ and $a_{i+1}wb_j$ are not in $L$ and $L'$.  This carries our induction hypothesis both  to  $i+1$ and $j$ and to $i$ and $j+1$. 

 In all cases our induction hypothesis is carried  to all $awb$, for $a$ and $b$ letters,
 and is now proved for all words of length $n+1$.
 \\

$(ii)$ implies $(iii)$, as a standard interval exchange transformation cannot have a wandering interval as it preserves the Lebesgue measure.\\

$(iii)$ implies the $\mathcal F$-flipped order condition by Proposition \ref{pord}. It implies recurrence as if a word $w$ of $L(T)$ is such that there is no $v\neq w$ beginning and ending with $w$, then the cylinder $[w]$ is a wandering interval.
\qed\\

  \begin{theorem}\label{tmin} Let $L$ be a language on an alphabet $\mathcal A$ and
$\mathcal F\subset \mathcal A$. Then $L$ is the language of a generalized (resp. standard) $\mathcal F$-flipped minimal interval exchange transformation if and only if it is aperiodic, uniformly recurrent, and satisfies an $\mathcal F$-flipped order condition. \end{theorem}
{\bf Proof}\\
For any language $L$, the minimality of the system
$(X_L, S)$ is equivalent to the uniform recurrence of $L$, see for example \cite{pyt} ch 5.

If $T$ is an $\mathcal F$-flipped interval exchange transformation, $L(T)$ satisfies the required order condition by Proposition \ref{pord}. If $T$ is minimal, $T$ has no periodic point, as finite orbits cannot be dense, thus $L(T)$ is aperiodic. The minimality of $T$ implies  the minimality of $X_{L(T]}$ as arbitrarily  long cylinders correspond to arbitrarily small intervals, thus $L(T)$ is uniformly recurrent.

In the other direction, we apply Theorem \ref{tg1} to get a standard interval exchange transformation for which $L=L(T)$. $(X_L, S)$ is minimal, thus $T$ can be not minimal only if arbitrarily long cylinders do not correspond to arbitrarily small intervals. This happens only if an infinite trajectory $u=u_0\cdots u_n\cdots $ corresponds to an interval $J$ in $(0,1)$. But either $T^ku=u$ for some $k>0$, which is excluded by aperiodicity, or else for all $k>0$ $T^kJ$ is disjoint from $J$ as they are in two disjoint cylinders;
but this is excluded as $T$ is a standard interval exchange transformation, preserving the Lebesgue measure. \qed\\

Both Theorem \ref{tg1} and Theorem \ref{tmin} should be compared with the characterizations in \cite{bc}. A particular case of Theorem \ref{tmin} is proved in \cite{fz3}, Theorem 2, for unflipped standard  interval exchange transformations. It uses a modified version of  the  {\em i.d.o.c. condition} introduced by Keane \cite{kea}; in the present context, this condition we use can be stated as: there is at least one point $\gamma_i$ (of Definition \ref{diet}), and there is no $i$, $j$, $k\geq 0$, such that $T^k\beta_i=\gamma_j$. Note that  this condition depends  on the defining intervals, and is not intrinsic to $T$ as the original 
Keane's condition.  The i.d.o.c. condition could be generalized  to flipped standard interval exchange transformations, but, to our knowledge, this has never been done, and it is not clear that any non-trivial example could satisfy it, as among flipped standard interval exchange transformations, the minimal ones constitute a small \cite{nog2} albeit non-empty \cite{nog1} set

Now  we can restate our old theorem in a cleaner form.  

\begin{theorem}\label{tido} A language  $L$ on at least two letters is the language of a standard unflipped interval exchange transformation satisfying the i.d.o.c. condition if and only if it satisfies an unflipped order condition, is aperiodic and uniformly recurrent, and has no connection. \end{theorem}

 This can be deduced from the present Theorem \ref{tmin} by noticing that our i.d.o.c.  condition implies Keane's one and thus minimality, and that for a minimal standard interval exchange transformation, our i.d.o.c. condition is just the  absence of connections in $L(T)$.\\

Note that minimal generalized interval exchange transformation cannot have wandering intervals; this can also be proved by using the above results: by minimality the existence of a wandering interval for $T$ implies the existence of  a wandering cylinder for the shift on $(X_{L(T)}, S)$, and this  contradicts the measure condition for $L(T)$.

\section{Generalized interval exchange transformations}\label{gai}

 We turn again to the study of languages satisfying an order condition or a weaker property: we are now interested in non-recurrent elements of $X_L$.

\begin{lemma}\label{nrec1}
 Suppose $L$ has no strong bispecial word. Let $x$ be in $X_L$. There exists  $m$ such that the infinite prefix $(x_n, n\leq -m)$ is left recurrent and the infinite suffix $(x_n, n\geq m)$ is right recurrent. \end{lemma}
{\bf Proof}\\
Suppose for all $m\geq m_0$ the infinite suffix $(x_n, n\geq m)$ is not right recurrent.
Then for such $m$ there exists a word $w^{(m)}$
which occurs exactly once in $(x_n, n\geq m)$.
As there is a finite number of letters, by going to a subsequence we can suppose  all the $w^{(m)}$ begin by the same letter, and thus this letter occurs infinitely often in $(x_n, n\geq m_0)$. Thus the longest prefix of  $w^{(m)}$ which occurs
 at least twice in $(x_n, n\geq m)$ is a nonempty $v^{(m)}$; then $v^{(m)}$ is right special, and, by applying Lemma \ref{orbls} and going to a subsequence, we can suppose $v^{(m)}$ is a suffix of a fixed infinite prefix, and also that the letter following $v^{(m)}$ in  $w^{(m)}$ is a fixed $a$. Now we can find $m<m'$ such that $v^{(m')}$ is not shorter than $v^{(m)}$, thus $v^{(m')}$ contains $v^{(m)}$, $v^{(m')}a$ contains $v^{(m)}a$, thus $v^{(m)}a$ occurs in $(x_n, n\geq m')$, which contradicts the assumption that $v^{(m)}$ is the longest possible. And similarly for the other assumption.\qed\\

\begin{proposition}\label{nrec2} Suppose $L$ has no strong bispecial  word. Then there are at most a finite number of {\em orbits} of $S$, namely sets $\{S^nx, x\in X_L\}$, such that $x$ is not  recurrent. \end{proposition}
{\bf Proof}\\
If $x$ is not recurrent, then some $(x_n, n\geq k)$ is not right recurrent.  Let $k$ be the largest of these $m$: it is finite by Lemma \ref{nrec1}, and necessarily the infinite suffix
$(x_n, n\geq k+1)$ has at least  two left extensions. Thus $(x_n, n\geq k+1)$  is one of the finitely many infinite suffixes $W_i$ of Lemma \ref{orbls}. If $(x_n, n\geq l)$ is one of the $W_i$ for a sequence $l\to -\infty$, then it is the same $W_i$ by going to a subsequence, thus $x$ is periodic and hence recurrent. Thus there is a smallest $m=m_0$ such that $(x_n, n\geq m)$ is one of the $W_i$, thus for $m<m_0$ there is only one left extension for $(x_n, n\geq m)$, and the sequence $x$ is known (up to shifts) if we know $W_i$. As there are finitely many $W_i$, this give finitely many orbits. \qed\\

\begin{lemma}\label{nrec3} Let $L$ be a non-recurrent language satisfying an $\mathcal F$-flipped order condition. We denote by $L'$, on  $\mathcal A'\subset \mathcal A$, the language generated by all the recurrent sequences in $X_L$. $L'$ is nonempty, recurrent, satisfies an $\mathcal F'$-flipped order condition for osme  $\mathcal F'\subset \mathcal F$,  and the orders defined on $X_L$ and $X_{L'}$ coincide on $X_{L'}$. \end{lemma}
{\bf Proof}\\
$L'$ is not empty
as there are recurrent sequences in $X_L$: indeed, by Lemma \ref{nrec1}, there are always right recurrent infinite suffixes $(x_n), n\geq k$ and these can be extended to the left as for each $m\geq k$ the word  $x_k... x_m$ occurs somewhere to the right, thus for some letter $x_{k-1}$ the word   $x_{k-1}x_k... x_m$ occurs in   $(x_n), n\geq k$  for infinitely many $m$, hence for all $m\geq k$. Thus we get a right recurrent infinite suffix $(x_n), n\geq k-1$, and by iterating the process we get a recurrent bi-infinite sequence $(x_n, n\in \GZ)$.

The other assertions are immediate as $L'$ is recurrent by construction and satisfies a restriction of the order condition on $L$. \qed\\

\begin{theorem}\label{tg2} A language $L$ is a natural coding of an $\mathcal F$-flipped generalized interval exchange transformation iff $L$ satisfies an $\mathcal F$-flipped order condition.\end{theorem}
{\bf Proof}\\ The ``only if" direction is Proposition \ref{pord}. We suppose now $L$ satisfies an  $\mathcal F$-flipped order condition. We suppose also that $L$ is not recurrent, otherwise we can conclude by Theorem \ref{tg1}. Let $L'$ be as in Lemma \ref{nrec3}.\\

Let $z$ be in $X_{L}\setminus X_{L'}$; by Lemma  \ref{lort} if the set $Z^-=\{t\in X_{L'}; t\leq z\}$ and $Z^+=\{t\in X_{L'}; t\geq z\}$ are nonempty, $Z^+$ has a unique
lower bound $z^+\in X_{L'}$ and $Z^-$ has a unique upper bound $z^-\in X_{L'}$; in that case, $z^+>z^-$, otherwise $z$ would be in $X_{L'}$.  If $\{t\in X_{L'}; t\leq z\}$ is empty, we define $z^-=z^+$ to be the  minimum of $X_{L'}$, and if  $\{t\in X_{L'}; y\geq z\}$ is empty, we define $z^+=z^-$ to be the maximum of $X_{L'}$.

We apply Theorem \ref{tg1} to find a standard interval exchange transformation $T'$, coded by $\mathcal A'$, such that $L'=L(T')$. $T'$ is not necessarily unique, we choose one. By Remark \ref{impp}, every sequence in $X_{L'}$ is either the trajectory of a point $x$ under $T'$, or an improper trajectory associated to a point $x$;  $x$ is not unique for periodic sequences, but for a given sequence the set of possible $x$ is an interval. For $z$ in  $X_{L'}\setminus X_L$, if  $z^-<z^+$,  there is no $t$ in $X_{L'}$ between them, this is possible only if $z^+$ and $z^-$ are two improper trajectories associated to the same point $x$; then $x$ is unique as it is the rightmost point of the interval associated to $z^-$ and the leftmost point of the interval associated to $z^+$, and we call it $p(z)$.   If $\{t\in X_{L'}; t\leq z\}$ is empty, we define $p(z)$ to be the leftmost point in the domain of $T'$, namely $0$, and if  $\{t\in X_{L'}; t\geq z\}$ is empty, we define $p(z)$  to be the rightmost point in the domain of $T'$, namely $1$.  In all cases, for a given $z$, for $|m|$ large enough $(S^mz)^-=(S^mz)^+$, which implies $p(S^{m+1}z)=T'p(S^mz)$.
 \\

By Proposition \ref{nrec2} and its proof, there are finitely many orbits in  $X_{L}\setminus X_{L'}$; we describe them as $S^nz_{(i)}, n\in \GZ, 1\leq i\leq K$, where $z_{(i)}$ is chosen so that $(z_{(i)m}, m \leq 0)$ is left recurrent but
$(z_{(i)m}, m\leq 1)$ is not left recurrent. We can do the same reasoning in both directions, thus $z_{(i)}=y_{(i)}w_{(i)}y'_{(i)}$ for an infinite prefix $y_{(i)}$ of a sequence in $X_{L'}$,  an infinite suffix $y'_{(i)}$ of a sequence in $X_{L'}$,  a (possibly empty) finite word $w_{(i)}$ in $L$; we know also that  $y_{(i)}=(z_{(i)m}, m\leq 0)$, but  $y'_{(i)}$ is not necessarily the first left recurrent infinite suffix of  $z_{(i)}$ as the latter might begin at a negative coordinate.

We look at the possible $S^nz_{(i)}$ such that $p(S^nz_{(i)})$ is a given point: this is a possibly infinite subset of $X_L$, we show that it has at most finitely many accumulation points for the usual topology: indeed, any accumulation point is a limit of points in the same orbit; if  $p(S^nz_{(i)})=p(S^{n'}z_{(i)})$ for $n\neq n'$, then  $y'_{(i)}$ is periodic, equal to $vvv...$,   $y_{(i)}$ is periodic, equal to $...v'v'v'$, and on this orbit there are at most finitely many  accumulation points, which are the iterates by $S$ of $...vvv...$ and $...v'v'v'...$

Thus for any given $x$ in $p(X_{L}\setminus X_{L'})$, all the $S^nz_{(i)}$ such that $p(S^nz_{(i)})=x$ can be ordered (by the order on $X_{L'}$) in a sequence $\zeta_n(x)$, $n\in H(x)$, $H(x)$ being  a finite or countable set. We now blow up  $x$ to a closed interval $\hat J(x)$ which is divided into adjacent open intervals $J(\zeta_n(x))$, ordered from left to right by the growing order on $\zeta_n(x)$, where the length of $J(\zeta_n(x))$ is $2^{-|m|}$ if $\zeta_n(x)=S^mz_{(i)}$.

Namely, for $x$  in $p(X_L\setminus X_{L'})$,  let $l(x)$ be the sum on $n\in H(x)$ of the lengths of $J(\zeta_n(x))$; $l(x)$ is  finite; for $x$ in $(0,1)$ and not in $p(X_{L}\setminus X_{L'})$, we put $l(x)=0$. For all $x$ in $(0,1)$, we define $\phi(x)= \sum_{y \in p(X_{L}\setminus X_{L'}), y<x}l(y)$, this last sum being also finite; let $M=1+\phi(1)+l(1)$. Then we send each point $x$ in $(0,1)$ to the interval (possibly reduced to a  point) $\hat J(x)=[x+\phi(x),x+\phi(x)+l(x)]$; for $l(x)\neq 0$ we divide $\hat J(x)$ into
$J_{\zeta_n(x)}$ as described above; the image of $(0,1)$ by this Denjoy-type  blow-up \cite{den} \cite{her} is $[0,M]$. On this interval we define $T$ by sending in an affine way $J(S^mz_{(i)})$ onto $J(S^{m+1}z_{(i)})$ for all $m$ and $i$, with a negative slope whenever $z_{(i)m}$ is in $\mathcal F$. If both $x$ and $T'x$ are not in $p(X_{L}\setminus X_{L'})$, we  define $T(x+\phi(x))=T'x +\phi(T'x)$. This leaves $T$ or $T^{-1}$ undefined on some points corresponding to $x$ where $T'$ or $T'^{-1}$ is undefined, on the endpoints of the $J(S^mz_{(i)})$ and finitely many images by $T'^n$ of these endpoints. but also  to $x$ (resp. $T'x$) not in $p(X_{L}\setminus X_{L'})$ but such that $T'x$ (resp. $x$) is in
$p(X_{L}\setminus X_{L'})$; the latter case will happen only for finitely many points, and then the  trajectory under $T'$ of either $x$ or $T'^x$ is some improper trajectory $(S^{m+1}z_{(i)})^+$ or $(S^{m+1}z_{(i)})^-$. Finally, we put in the interval $I_e$, $e\in \mathcal A$, all points $x+\phi(x)$ for $x$ not in $p(X_{L}\setminus X_{L'})$ such that $x$ is in the defining interval $I'_e$ of $T'$, and all the $J(S^mz_{(i)})$ for which $e=z_{(i)m}$. By the above considerations, all trajectories for $T'$ remain as actual or  improper trajectories for $T$.

Because we have respected the order, $T$ is indeed a generalized interval exchange transformation with the required orders and $\mathcal F$, and its language is indeed $L$, as the sequences in $X_L$ are the actual or improper trajectories for $T$. \qed\\

\begin{remark}\label{2bu} The blow-up we made in the above proof sends points $x$ to intervals $[x+\phi(x),x+\phi(x)+l(x)]$, while \cite{cage} \cite{bhm} \cite{co} \cite{mmy} use intervals $[\phi(x),\phi(x)+l(x)]$ instead; the latter kind of blow-up is needed to build affine $T$, but two different points $x$ and $x'$ may be sent to the same point, though this will not happen if $T'$ is minimal as then between two different points there is always some $p(z)$, while in the former case the blow-up is always injective but does not produce affine interval exchange transformations. In the proof of Theorem \ref{paf} below we shall use a mixture of these two.\end{remark}

We shall see how this proof works on Examples \ref{fs2} of Section \ref{ce2} and  \ref{mon} of Section \ref{ce3}.

\section{Affine interval exchange transformation}\label{sai}

\begin{theorem}\label{paf} If $L$ is a natural coding of an $\mathcal F$-flipped affine interval exchange transformation for which the absolute value of the slope is $\exp{\theta_e}$ on the defining interval $I_e$, then $L$ satisfies an  ${\mathcal F}$-flipped order condition and for each non recurrent sequence $z$ in $X_L$, $\sum_{n\geq 0}\exp{\sum_{j=0}^n \theta_{z_j}}<+\infty$, and  $\sum_{n>0}\exp{-\sum_{j=-n}^{-1} \theta_{z_j}}<+\infty$.

 If $L$ satisfies an ${\mathcal F}$-flipped order condition and  there exist real numbers $\theta_e, e\in\mathcal A$, such that for each non recurrent sequence $z$ in $L$, $\sum_{n\geq 0}\exp{\sum_{j=0}^n\theta_{z_j}}<+\infty$, and  $\sum_{n>0}\exp{-\sum_{j=-n}^{-1} \theta_{z_j}}<+\infty$, then $L$ is a group coding of an ${\mathcal F}$-flipped  affine interval exchange transformation.  \end{theorem}
{\bf Proof}\\
Let $z$ be a non recurrent sequence. Let $L'$ be as in the proof of of Theorem \ref{tg2}. As  $z$ is not recurrent, $z$ cannot have all its words in the same  uniformly recurrent component (of   Proposition \ref{cmin}) of $L'$, and thus not  all its words  in $L'$. Thus we can find a word $w=w_1...w_k$ which is not in $L'$ and appears only once in $z$; by Lemme \ref{nrec1} it will appear only finitely many times in each of the other non recurrent sequences, and by Lemma \ref{nrec2} we can extend it so that it appears only in $z$ and its shifts. Thus the cylinder $[w]$ is disjoint from all its iterates by $T$, and  for the Lebesgue measure $\mu$ both  $\sum_{n\geq 0}\mu(T^n[w])$ and $\sum_{n\geq 0}\mu(T^{-n}[w])$ must be finite. Suppose $w=z_k...z_{k+k'}$; then $T$ is of slope $\exp{\theta_{z_k}}$ on $[w]$, $\exp{\theta_{z_{k+1}}}$ on $[Tw]$, $\exp{\theta_{z_{k-1}}}$ on $T^{-1}[w]$ and so on, thus $\mu(T[w])=\exp{\theta_{z_k}}\mu[w]$, $\mu(T^2[w])=\exp{\theta_{z_{k+1}}}\mu(T[w])$,  $\mu([w])=\exp{\theta_{z_{k-1}}}\mu(T^{-1}[w])$ and so on. Our conditions come from replacing $k$ by $0$, which amounts to multiplying by a constant.\\

We start now from $L$ and build a generalized interval exchange transformation $T$, as in the proof of Theorem \ref{tg2}, but with two modifications.

First, the length we assign to the interval $J(S^mz_{(i)})$ is $1$ for $m=0$, then is defined such that $J(S^mz_{(i)})$ is sent to $J(S^{m+1}z_{(i)})$ by an  affine map of slope $\exp \theta_{z_{(i)m}}$ if  $z_{(i)m}$ is in $\mathcal F^c$, $-\exp \theta_{z_{(i)m}}$ if  $z_{(i)m}$ is in $\mathcal F$.

Second, when we make the blow-up of $(0,1)$ to the domain of $T$, we replace the map $x\to \phi(x)+x$ by the map $x\to \phi(x)+1_K(x)x$ for some union of intervals $K$. To find $K$, we apply the proof of Proposition \ref{cmin} to the recurrent language $L'$ : for some fixed $N$ large enough, its Rauzy graph $G_N$ splits into a finite number of connected components, corrresponding to languages which  either come from one purely periodic orbit or are aperiodic and uniformly recurrent. This implies that  there exist a finite number of $T'$-invariant sets $K_j$, which are finite unions of cylinders of length $N$, and  such that, on each $K_j$, either $T'$ is minimal, or every point is periodic (note that this reproves  Corollary 2.9 of \cite{arth}, and  extends it to flipped interval exchange transformations). 

Now, for each $i$, we look at those points   $p(S^mz_{(i)}), |m|\geq m_i$, where $m_i$ is chosen such that for these points  $p(S^{m+1}z_{(i)})= T'p(S^mz_{(i)})$, and let $K^c$ be the union of those $K_j$ which contain at least one  $p(S^{m_i}z_{(i)})$ or  $p(S^{-m_i}z_{(i)})$ in their closure.  If $x\neq x'$, our construction ensures that $x$ and $x'$ are blown-up into different points, unless maybe if $x$ and $x'$ belong to the same connected component $K'$ of the same $K_i$: in the latter case either minimality implies that there is at least one $J(S^mz_{(i)})$ between their images, in which case again $x$ and $x'$ are blown-up into different points, or periodicity implies that $x$ and $x'$ have the same trajectory, and this trajectory is not lost in the blow-up as it remains as the improper trajectory of some endpoint of $K'$. Thus again all trajectories for $T'$ remain as actual or  improper trajectories for $T$.

The map $T$ built in this way is a generalized $\mathcal F$-flipped interval exchange transformation with defining intervals $I_e$. By the same proof as in \cite{cage}  \cite{bhm} \cite{co} \cite{mmy}, if the set $K$ in $(0,1)$ is blown up to $\hat K$, $T$ is affine of slope  $\exp{\theta_e}$ (resp. $-\exp{\theta_e}$ if $e$ is in $\mathcal F$) on each connected component of $I_e\cup \hat K^c$, while  each connected component of $I_e\cup \hat K$ is the union of finitely many $J(S^mz_{(i)})$, on which $T$ is affine of slope  $\exp{\theta_e}$ (resp. $-\exp{\theta_e}$ if $e$ is in $\mathcal F$) and finitely many intervals on which $T$ is affine of slope $1$ (resp. $-1$ if $e$ is in $\mathcal F$). Thus by dividing the intervals further we can make $T$ affine with constant slope  on its defining intervals, and $L$ is a grouped coding of an affine interval exchange transformation.  \qed\\

\begin{corollary}\label{caf}
 If $L$ satisfies an ${\mathcal F}$-flipped order condition and $X_L$ contains some non-recurrent sequences,  a sufficient condition
  for $L$ to be  a natural  coding of an ${\mathcal F}$-flipped  affine interval exchange
  transformation is that 
  \begin{itemize} 
  \item there exist real numbers $\theta_e, e\in\mathcal A$, such that for each non recurrent sequence $z$ in $L$, $\sum_{n\geq 0}\exp{\sum_{j=0}^n \theta_{z_j}}<+\infty$, and  $\sum_{n>0}\exp{-\sum_{j=-n}^{-1} \theta_{z_j}}<+\infty$, 
  \item each factor of a recurrent sequence in $X_L$ is also a factor of a non-recurrent sequence in $X_L$.
  \end{itemize} 
\end{corollary}
{\bf Proof}\\
With the notations of the proof of  Theorem \ref{paf}, the last condition implies that in the closure of  every $K_j$ there is at least one  $p(S^{m_i}z_{(i)})$ or  $p(S^{-m_i}z_{(i)})$, thus $K$ is empty and there is no subinterval of $I_e$ on which $T$ is affine of slope $1\neq \exp{\theta_e}$ or $-1\neq -\exp{\theta_e}$, which ensures that $L$ is the natural coding of $T$. \qed\\

This criterion will be further discussed in Section \ref{cev} below.\\

We can get a  criterion for $L$ to be a grouped coding of an affine interval exchange transformation by  essentially  rephrasing  Theorem \ref{paf}.

\begin{definition}\label{lgc} Suppose $L$ satisfies an $\mathcal F$-flipped order condition.  A language $\hat L$ is a {\em splitting} of $L$ if its letters are  $e_i$, $e\in\mathcal A$, $1\leq i\leq k_e$, forming the set $\hat{\mathcal A}$, such that
$L=\phi(\hat L)$ where we define $\phi(e_i)=e$ for all $e$ and $i$,  and extend $\phi$ to a morphism for the concatenation on words, and
$\hat L$ satisfies an $\hat{\mathcal F}$-flipped order condition, where
\begin{itemize} \item $\hat{\mathcal F}$ is made of all the $e_i$ such that $e$ is in  $\mathcal F$,
  \item if $e$ is in $\mathcal F^c$, $e_1$, ...  $e_{k_e}$ are consecutive and  ordered in the same way by the orders $<_A$  and  $<_D$ of $\hat L$,
  \item if $e$ is in $\mathcal F$, $e_1$, ... , $e_{k_e}$ are consecutive and ordered in opposite ways by the orders $<_A$  and  $<_D$ of $\hat L$,
  \item if $e\leq_A e'$, resp.  $e\leq_D e'$  in $L$, then  $e_i\leq_A e'_j$, resp.  $e_i\leq_D e'_j$,  in $L'$ for all $i$ and $j$,
  \end{itemize}
\end{definition}

\begin{theorem}\label{tg3} $L$ is a group coding of an ${\mathcal F}$-flipped  affine interval exchange transformation if and only  $L$ satisfies an  ${\mathcal F}$-flipped order condition and  there exists a splitting $\hat L$ of $L$,  and real numbers $\theta_e, e\in\hat{\mathcal A}$, such that for each non recurrent sequence $z$ in $X_{\hat L}$, $\sum_{n\geq 0}\exp{\sum_{j=0}^n \theta_{z_j}}<+\infty$, and  $\sum_{n>0}\exp{-\sum_{j=-n}^{-1} \theta_{z_j}}<+\infty$.  \end{theorem}

{\bf Proof}\\ In one direction, we suppose $L$ is a grouped coding of $T$  affine,  define $\hat L$ to be the natural coding of $T$ and apply  Theorem \ref{paf}. In the other direction, we start from $\hat L$ and use Theorem \ref{paf} to build an affine interval exchange transformation $T$ such that $\hat L$ is a grouped coding of $T$. But then $L=\phi(\hat L)$  is a further grouped coding of $T$. \qed\\

\section{Counter-examples and questions}\label{sexa}
\subsection{Affine with natural coding}\label{ce1}

 \begin{example}\label{fs1} Let $L$ be generated  by the bi-infinite sequence $\ldots 111222 \ldots $. Note that it is of complexity $n+1$ but not uniformly recurrent, and in the founding paper \cite{hm} it is not included in Sturmian languges, we can call it a {\em fake Sturmian} language.  It satisfies the unflipped order condition with $1<_D2$, $1<_A2$, but (unsurprisingly as it is not recurrent)  is not the language of a standard interval exchange transformation as that could only be the identity on two disjoint open intervals $I_1$ and $I_2$, and the only possible words are $1^n$ and $2^m$. However, $L$ is the natural coding of an affine interval exchange transformation: $L_2$ is the language of length $2$ of any affine $2$-interval exchange transformation, with the same orders, such that $TI_1$ is strictly longer than $I_1$, and, as $L$ is determined by $L_2$ because  there is no bispecial word except the empty one, $L$ is indeed the  natural coding of any of these affine interval exchange transformations.  \end{example}

  Example \ref{fs1} can also be dealt with as in Theorem \ref{paf}, with $T'$ being the identity acting on $I_1$ and $I_2$.\\
	
	Because of Theorem \ref{paf} and Corollary \ref{caf}, if we want  to build nontrivial examples of affine interval exchange transformation with a non recurrent language $L$ as a natural coding, we can just use \cite{cage} \cite{bhm} \cite{co} \cite{mmy}: starting from a minimal standard interval exchange transformation $T'$, we find  a point $x_{\star}$ and a map $\theta$ constant on the defining intervals  such that the Birkhoff sums
$\sum_{j=0}^n \theta (T'^jx_{\star})$  and $\sum_{j=-1}^{-n} \theta (T'^jx_{\star})$ behave as we need, then  make a blow-up of all points in the orbit of $x_{\star}$, and get an affine $T$ with a non-recurrent  trajectory using an extra letter as in Example \ref{mon}; note however that, because of the order condition,  $x_{\star}$ must be  an endpoint of a defining interval of $T'$; if  $x_{\star}$ is not already such an endpoint, we can make it be one with a modification of $T'$, by splitting the interval $I'_e$ containing $x_{\star}$ into $I'_e\cap (0, x_{\star})$ and $I'_e \cap (x_{\star},1)$.  We can also get examples by  blowing up two half-orbits, of two points $x_{\star}$ and $x'_{\star}$ so that the Birkhoff sums $\sum_{j=0}^n \theta (T'^jx_{\star})$ and $\sum_{j=-n}^{-1} \theta (T'^jx'_{\star})$
behave as we need, provided we then modify $T'$ so that $x_{\star}$ and $x'_{\star}$ are endpoints of defining intervals; this is always possible if $T'$ is minimal but non uniquely ergodic, as we can find two invariant measures $\nu$ and $\nu'$ and a vector $\theta$ so that $\sum\theta_e\nu(I_e)<0$, $\sum\theta_e\nu'(I_e)>0$, and then choose $x_{\star}$ generic for $\nu$, $x_{\star}$ generic for $\nu'$. But if $T'$ is uniquely ergodic, then we need $\sum\theta_e\nu(I_e)=0$ for the invariant measure $\nu$, and, by \cite{hal}, for any given $\theta$ the set of $x$ such that $\sum_{j=0}^n \theta (T'^j)\to -\infty$  or $\sum_{j=-1}^{-n} \theta (T'^jx)\to+\infty$ is of measure $0$.
  
\subsection{Affine with grouped coding}\label{ce2}
 \begin{example}\label{fs2} Let $L$ the language  generated  by the bi-infinite sequence $\ldots 1112111 \ldots $, which is a {\em skew Sturmian} language as defined in \cite{hm}.  It satisfies the unflipped order condition with $1<_D2$, $2<_A1$, but  is not the natural coding of any affine interval exchange transformation $T$: indeed, the sequence $\ldots 1111 \ldots$ in $X_L$ would define a fixed point $x$ for $T$, in the interior of $[1]$, and, if $0<y<x$ is the right endpoint of $T[2]$, $T$ would have to send $(0,x)$ to  $(y,x)$ and $(x,1-y)$ to  $(x,1)$, thus having a slope $<1$ on  a part of $[1]$ and a slope $>1$ on another part. However, if $\tilde L$ is the language generated  by the bi-infinite sequence $\ldots 3332111 \ldots $, as in Example \ref{fs1} $\tilde L$ is the natural coding of any affine interval exchange transformation $T$ sending $I_1=(0,x)$ to  $(y,x)$, $I_3=(x,1-y)$ to  $(x,1)$, $I_2=(1-y,1)$ to  $(0,y)$, with $0<y<x<1-y$. If we now code $T$ by the intervals $\tilde I_1=I_1\cup I_3$ and $\tilde I_2=I_2$, we see that $L$ is indeed a grouped coding of an affine interval exchange transformation as in Definition \ref{gc}. \end{example}

We show how the proof of Theorem \ref{tg2} works  on Example \ref{fs2}. $L'$ is the language generated by the single sequence $y=...111...$, and thus $T'$ is the identity acting on $I_1=(0,1)$. In the single non recurrent orbit $...11211...$, $z$ is the point for which $z_0=2$, and we have $Sz<S^2z<S^3z<...<y<z<...
 S^{-3}z<S^{-2}z<S^{-1}z<z$, with $y$ being the limit of the $S^nz$ when $n\to\pm\infty$. Thus $p(S^nz)=1$ for all $n\leq 0$, $p(S^nz)=0$ for all $n\geq 1$; our construction gives an affine interval exchange transformation $T$ sending $(0,1)$ to $(1/2,1)$, $(1,2)$ to itself, $(2,3)$ to $(2,4)$, $(3,4)$ to $(0,1/2)$, while the cylinder $[1]$ is the interval $(0,3)$ and the cylinder $[2]$ is the interval $(3,4)$. $L$ is the natural coding
 of $T$ if we consider it as a generalized interval exchange transformation, but is only a grouped coding of an affine interval exchange transformation, as remarked above. Note that the fact that $T$ is affine is a (conceivably) accidental result of our choices of the lengths of the added intervals, see the proof of Theorem \ref{paf}, and that the present map $T$ is slightly different from the one of Example \ref{fs2}; the latter could be obtained from the former by collapsing the interval $(1,2)$, as in the second blow-up of Remark \ref{2bu} (this  is possible though $T'$ is not minimal because  all the points in $(1,2)$  have the same trajectory, $...1111...$, and this will remain,  as an improper trajectory, after the collapsing).

 \subsection{Natural versus grouped}\label{cev}
 \begin{conj}\label{con1} The conditions in Theorem \ref{paf} are necessary and sufficient for $L$ to be a natural coding of an ${\mathcal F}$-flipped affine interval exchange transformation. \end{conj}

Indeed, the last condition of Corollary \ref{caf} is satisfied if $T'$ is minimal, but also by Example \ref{fs1} above. But
 this condition is {\it not necessary}, for different reasons shown in the following three examples, most significantly Example \ref{exaf}. 
  \begin{example} Let $L$ be the language generated by the two infinite sequences $...111222....$ and $...3333...$. Then  $K$ is a union of defining intervals of $T'$, so the grouped coding of Theorem \ref{paf} turns out to be a natural coding. \end{example} 
  \begin{example} Let  $L$ be generated by $...121312131213...$ and $..114114114115115115...$ Then 
 $K$ is nonempty and not a union of defining intervals, but an allowed value of
$\theta_1$ is $0$, and $1$ is the only $e$ for which   $I'_e$  contains strictly a connected component of $K$, thus again the proof of Theorem \ref{paf} provides a natural coding. \end{example}
\begin{example}\label{exaf} Let $L$ be generated by $...332332331331....$ and $...413241324132...$; then either $\theta_1\neq 0$ or $\theta_2\neq 0$, hence the proof of Theorem \ref{paf} would yield an interval exchange transformation with two different slopes, $1$ and $\exp{\theta_e}$, on $I_e$ for $e=1$ or $e=2$, thus $L$ is not the natural coding of any affine interval exchange  transformation built by the proof of Theorem \ref{paf}. However,  $L$ is the natural coding of  another affine interval exchange transformation, namely any member of the family defined by the orders $1 <_D 2<_D 3<_D 4$ and $4 <_A 3 <_A 1 <_A 2$, slopes $a$ on $I_1$, $1/a$ on $I_2$, $1$ on $I_3$ and $I_4$, and lengths $|I_1|=l+r$, $|I_2|=a(l+r)$, $|I_3|=al+2(a+1)r$, $|I_4|=l$, for any given $0<a<1$, $l>0$ and $r>0$.\end{example} 

The only example we know of a grouped coding which is not a natural coding is Example \ref{fs2}, which does not satisfy the condition on Birkhoff sums: whatever the value  of $\theta_1$ and $\theta_2$, either the $\exp n\theta_1$, $n\geq 0$, or the $\exp n\theta_1$, $n<0$, must be bounded away from $0$ (while for the language $\tilde L$ of Example \ref{fs2} we can take $\theta_1<0$ and $\theta_3>0$).

\begin{question}\label{qaprd} Does there exist an aperiodic language which is a grouped coding of an affine interval exchange transformation, but not a natural coding of any affine interval exchange transformation? \end{question}

Questions \ref{con1} and \ref{qaprd}  suggest what we dare not call a conjecture.

\begin{question}\label{con2} Is it true that $L$ is a group coding of an ${\mathcal F}$-flipped  affine interval exchange transformation if and only if $L$ satisfies an ${\mathcal F}$-flipped order condition and  there exist real numbers $\theta_e, e\in\mathcal A$, such that
the two following conditions hold?
\begin{itemize}\item  For each non recurrent sequence $z$ in $L$ which is not ultimately periodic to the left, \\  $\sum_{n\geq 0}\exp{\sum_{j=0}^n\theta_{z_j}}<+\infty$.
 \item For each non recurrent sequence $z$ in $L$ which is not ultimately periodic to the right, \\  $\sum_{n>0}\exp{-\sum_{j=-n}^{-1} \theta_{z_j}}<+\infty$. \end{itemize} \end{question}

\begin{remark} As in the last example above, the same language  can be  a coding of two fairly different affine interval exchange transformations. In particular, Sturmian languages are languages of standard interval exchange transformations, thus also of affine interval exchange transformations without wandering intervals. It is known since \cite{hm} that if we want to generate a given Sturmian language $L$ by a standard $2$-interval exchange transformation, the parameter $\mu[I_2]$ has two possible values, or one up to reversal of the orientation. If we want the interval exchange transformation to be affine, we have two parameters, for example its slopes $a$ and $b$ on $I_1$ and $I_2$, and we do not know which amount of freedom we have in this case.\end{remark}

  \begin{question}\label{sta} For a given $L$, what can be said of  the set of $(a,b)$ for which $L$ is the language of the corresponding affine interval exchange transformation? \end{question}

  About this question,  it is worth mentioning a surprising example due to M. Shannon (unpublished):

  \begin{example}\label{sha} One can build (by geometrical methods) a grouped coding of an affine interval exchange transformation $T$ which is a Sturmian language, associated to a rotation of angle $\alpha$, but while $\alpha$ is in $\GQ(\sqrt3)$, all the parameters defining $T$ are in $\GQ(\sqrt5)$. \end{example}

\subsection{Generalized}\label{ce3}
 \begin{example}\label{mon} Let $L'$ be  the Sturmian language which is  the natural coding of the unflipped standard  interval exchange transformation $T'$
sending $I_1=(0,1-\alpha)$  to $(\alpha,1)$ and $I_2=(1-\alpha, 1)$ to $(0,\alpha)$ for an irrational $\alpha<1/2$. Let $y_n=i$ whenever $T^n\alpha$ is in $I_i$, $n\geq 0$, and  $y'_n=i$ whenever $T^n(1-2\alpha)$ is in $I_i$, $n\leq 0$; when $\alpha=\frac{3-\sqrt 5}2$, $y$ is the so-called Fibonacci sequence on $1$ and $2$, and $y'$ is $y$ written backwards. Let $L$ be the language generated by the infinite sequence $...y'_{-2}y'_{-1}y'_03y_0y_1y_2...$  It satisfies the unflipped order condition with $1<_D3<_D2$, $2<_A3<_A1$ (note that no other unflipped order is possible, because of the way the empty bispecial is resolved, nor is any flipped order because of the way the bispecial $0$ is resolved).

 By Theorem \ref{tg2}  $L$ is the natural coding of a generalized interval exchange transformation, but it is not the natural or grouped coding of any affine interval exchange transformation: this will be a straightforward consequence of either one of  two independent results we show below, Theorems \ref{pmon} and  Theorem \ref{ttrok}. \end{example}

 Example \ref{mon} is a typical case to see how the proof of Theorem \ref{tg2} works. Then $L'$ and $T'$ are as described in its definition. Then the only not right  recurrent orbit in $X_L$ is the one of $y'3y$ defined above. There is one $z$, the element of this orbit such that $z_0=3$. Because of the order condition, the defining  interval $I_3$ of $T$ must be between $I_1$ and $I_2$, and $TI_3$ between $TI_2$ and $TI_1$; this implies that $p(z)=1-\alpha$ and $p(Tz)=\alpha$. And indeed
   $z^+$ is the element of the orbit $y'12y$ such that $z^+_0$ is the $1$ just after $y'$,
   $z^-$  is the element of the orbit $y'21y$ such that $z^-_{0}$ is the $2$ just after $y'$;
similarly
 $(Tz)^+$  is the element of the orbit $y'21y$ such that $(Tz)^+_{-1}$
  is the $1$ just before $y$, $(Tz)^-$  is the element of the orbit $y'12y$ such that $(Tz)^-_{-1}$ is the $2$ just before $y$. As there is no other discontinuity $p(T^nz)=T'^{n-1}\alpha$ for $n\geq 1$, $p(T^{-n}z)=T'^{-n}(1-\alpha)$ for $n\geq 1$. Thus we make the above construction by replacing each point $T^m\alpha$, $m\geq 0$,  by an interval $\hat J(T^m\alpha)$  of length $2^{-m-1}$, the point $1-\alpha$ by an interval $\hat J(1-\alpha)$  of length $1$, each point $T^{-m}(1-2\alpha)$, $m\geq 0$, by an interval $\hat J(T^{-m}(1-2\alpha))$  of length $2^{-m-1}$. \\

\begin{theorem}\label{pmon} Let $L$ be non recurrent, and
a natural coding  of an unflipped  generalized interval exchange transformation $T$. Suppose the language $L'$ of Lemma \ref{nrec3} is aperiodic, uniformly recurrent, and  its arrival and departure order are conjugate by a circular permutation. Then $T$ cannot be {\em of class P}, class P \cite{her} meaning that, except on a countable set of points, $DT$ exists and $DT=h$ where $h$ is a function with bounded variation, and $\ab h     \ab$ is  bounded from below by a strictly positive number. \end{theorem}
{\bf Proof}\\
A sequence in $X_L$  can be a trajectory for $T$ (actual or improper) either of a single point or of all points in an interval; the intervals which correspond to a single trajectory include those coded by non-recurrent trajectories,  but  possibly  others; as intervals corresponding to different trajectories are all disjoint, there are at most countably many of them; let $E_n$, $n\in \GZ$, be these intervals.  We build another interval exchange transformation $\hat T$ on the interval $(0,1)$ by making a deblow-up: this will be done in two different ways corresponding to the inverses of  the two blow-ups in Remark \ref{2bu}. If $\sum_{n\in\GZ}\mu(E_n)<1$, $\mu$ being the Lebsegue measure, we just shrink each $E_n$ (including its endpoints) to a single point, by sending $x$ to $x-\mu(]0,x[\cup\cap_{n\in\GZ}E_n)$, sending $(0,1)$ to a smaller interval, which we rescale to get $(0,1)$. If $\sum_{n\in\GZ}\mu(E_n)=1$, we use the semi-conjugacy of \cite{cage} \cite{bhm} \cite{co} \cite{mmy}: let $x_n$ be the left end of $E_n$; the condition on the lengths of $E_n$ implies the $x_n$ are dense in $(0,1)$; then we send to $x_n$ all points in $E_n$, while a point $x$ which is not im any $E_n$ is sent to the
unique (by density of the $x_n$) point $x'$ such that
$\sum_{x_n<x'}\mu(E_n)=x$.   In both cases all intervals coded by non-recurrent  trajectories are shrunk to single points, and $\hat T$ has no wandering interval. Then the natural coding of $\hat T$ is the language $L'$, as every recurrent trajectory of $T$ remains as an actual or improper trajectory of $\hat T$, while non-recurrent trajectories of $T$ are not trajectories of $\hat T$.

 By construction, each trajectory in $X_{L'}$ is the (actual or improper) trajectory under $\hat T$ of one point in $(0,1)$, and also the (actual or improper) trajectory under $T'$ of one point in $(0,1)$. By Remark \ref{ordint} the order on  $X_{L'}$, whether we consider trajectories of $T'$ or $\hat T$,  corresponds to the natural order on $(0,1)$, thus we can make a bijective, continuous and increasing map $\psi$  from $(0,1)$ to $(0,1)$ such that  $\psi\circ \hat T=T'\circ \psi$.\\

 Suppose $T$ is  of class P. Let $J$ be a wandering interval for $T$; then 
 we have $\mu(T^nJ)=\mu(J)\exp{\sum_{j=0}^{n-1} \theta_j}$ where $\exp\theta_j$ is the average value of $DT$ on $T^jJ$. If $h$ is the function defining class $P$, we define a function $\theta$,  by $\theta(x)=\log h(x)$ when $x$ is not in any $E_n$, while, for any $x$ in $E_n$, $\theta(x)$ is the logarithm of the average value of $h$ on $E_n$. Then the measures of $T^nJ$ are given by
  for $n>0$, $\mu(T^nJ)=\mu(J)\exp{\sum_{j=0}^{n-1} \theta(T^jx)}$ for some point $x$ in $J$, while  for $n<0$, $\mu(T^nJ)=\mu(J)\exp{-\sum_{j=n}^{-1} \theta(T^jx'})$ for some point $x'$ in $J$. As $J$ is a wandering interval, $\mu(T^nJ)$ tends to $0$ when $n$ goes to $\pm\infty$.
  
   The deblow-up sends Birkhoff sums of $\theta$ under $T$ to Birkhoff sums of a  function $\hat \theta$ under $\hat T$,  thus by the hypothesis
 these are Birkhoff sums of a measurable function under an irrational rotation, and, as in the case of a natural coding, because of the  unique ergodicity, we must have $\int\hat\theta d\mu=0$ for the Lebesgue measure.  The properties of $h$ imply that the function $\theta$ has bounded variation, thus so does $\hat\theta$; but then  by the Denjoy-Koksma inequality \cite{her} the Birkhoff sums  are bounded on a subsequence corresponding to the denominators of the partial quotients of the angle, and we  get a contradiction.\qed\\
 
The above theorem applies to Example \ref{mon} or any example built in the same way from any Sturmian word, but also, in contrast with Theorem \ref{ttrok} or Corollary \ref{csau} below, we may equip $T'$ with extra points $\gamma_i$ which are not discontinuities, giving extra possibilities for building $L$. All these give counter-examples of languages which are natural codings of a generalized interval exchange transformation, but not natural codings of any generalized interval exchange transformation of class P, and thus not  grouped codings of any affine interval exchange transformation. \\

\begin{theorem}\label{ttrok} Let $L'$ be a natural coding of a non purely periodic unflipped standard  interval exchange transformation. Let $w_n=av_nb$, $a\in \mathcal A$, $b\in \mathcal A$, be an infinite sequnece of bispecial words in $L'$ . Let $u$ be the infinite prefix in $X_{L'}$ ending with $w_n$ for all $n$, and $v$ the infinite suffix beginning with $w_n$ for all $n$. Let $\omega$ be a symbol which is not a letter of $L'$, and  $L$ be the language generated by the union of all words in $L'$ and the bi-infinite word $u\omega v$.
 
 Then $L$ is a natural coding of a generalized interval exchange transformation, but not 
 a  grouped coding of any affine interval exchange transformation. \end{theorem}
{\bf Proof}\\
We build a generalized interval exchange transformation with natural coding $L$ as in the proof of Theorem \ref{tg2}, by adding a point coded by $\omega$ at position $\gamma$, where $T'^{-|w_n|}\gamma$ is in $[w_n]$ for all $n$,  and its image at position $\beta$, where $\beta$ is in $[w_n]$ for all $n$, and making a Denjoy-Koksma blow-up.\\
 
 Suppose $L$ is a  grouped coding of an affine interval exchange transformation $T$.Let $\exp\theta(x)$ be the slope of $T$ at point $x$, $\theta$ is piecewise constant with  $K\geq 0$ jumps. As $[\omega]$ is a wandering interval, we can write the $T^j[\omega]$ as a (Rokhlin) tower, in which each discontinuity of $\theta$ appears at most once. Thus $[\omega]$ is partitioned into at most $K+1$ disjoint intervals, on which for any $n>0$ $\sum_{j=-n}^{-1}\theta(T^jx)$ and $\sum_{j=0}^{n}\theta(T^jx)$ are constant; we choose one of these subintervals, which is also a wandering interval, and call $Q_{-n}$ and $Q_n$ the respective value of these sums on it; when $n\to +\infty$, $Q_{-n} \to +\infty$ and $Q_n \to -\infty$ .

We look now at the induced (or first return) map of $T$ on the interval $[w_n]$, for fixed $n$. By the standard reasoning of \cite{kea}, this interval is partitioned  into $R$ subintervals $J_{i,n}$, which are the bases of $R$ disjoint  Rokhlin towers made with disjoint intervals $J_{i,n}$, ... $T^{h_{i,n}-1}J_{i,n}$, and then 
$T^{h_{i,n}}J_{i,n}$ is in $[w_n]$. The levels $T^lJ_{i,n}$ are in the defining interval corresponding to the $l+1$-th letter of the word $w_n$, and their union for fixed $l$ is an interval for $l\leq h'_n-1$, where  $h'_n=\min_i h_{i,n}$. Among these towers, there is one for which $h_{i_0,n}=|w_n|+1$
and $T^{h_{i_0,n}-1}[w_n]$ is in $[\omega]$. All the others are coded by orbits of $T'$, and $h'_n=\to+\infty$ when $n\to +\infty$, otherwise the $w_n$ would be in a periodic orbit, which is impossible as they are bispecial. \\

As $\omega$ is always followed by $w_n$, there is a subinterval of $[w_n]$ on which 
$\sum_{j=-1}^{l-1}\theta(T^jx) = Q_l$ for all $l\leq h'_n$. 

For fixed $n$, each discontinuity of $\theta$ appears at most once in the union of towers of bases $J_{i,n}$, and all these towers are adjacent from levels $0$ to $h'_n-1$; thus, for two different points $x$ and $x'$ in $[w_n]$  the values of $\theta(T^jx)$ and $\theta(T^jx')$ are different for at most $K$ distinct values of $j$ between $0$ and $h'_n-1$. Thus, taking also into account the replacement of $T^{-1}x$ by $T^lx$,  there exists a constant $M$ such that for all $n$,  all $l<h'_n$, all $x$ in $[w_n]$, $Q_l-M < \sum_{j=0}^{l}\theta(T^jx)< Q_l+M$. 

As $\omega$ is always preceded by $w_n$, there is a subinterval of $J_{i_0,n}$ on which 
$\sum_{j=0}^{h_{i_0,n}-2}\theta(T^jx) = Q_{-h_{i_0,n}+1}$. This implies that $Q_{-h_{i_0,n}+1}$ is between $Q_{h_{i_0,n}-2}-M$ and $Q_{h_{i_0,n}-2}+M$ for all $n$, which is a contradiction because of their behaviour when $n\to+\infty$.\qed\\

This gives us many counter-examples, of natural codings of a generalized interval exchange transformation which are  not grouped codings of any affine interval exchange transformation, including Example \ref{mon} for completely different reasons as the ones in Theorem \ref{pmon}, but for a given $T'$ a priori we cannot use all possible positions to add an interval $[\omega]$, nor can we  create new points $\gamma_i$ as possible positions.

\subsection{Geometric}\label{pascal}
In an intermediate case between the rotations of Theorem \ref{pmon} and the general systems of Theorem \ref{ttrok}, we give counter-examples which are not natural codings of an affine interval exchange transformation, starting from some famous $T'$, in which we can use all possible positions $\gamma_i$ to add an interval, but not create new points $\gamma_i$ as possible positions. 
\medskip

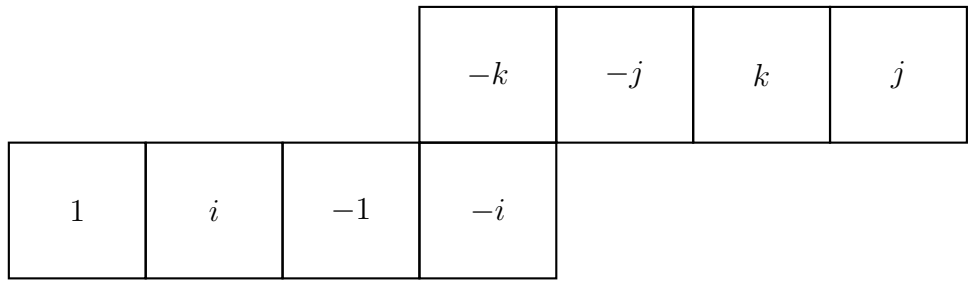
\begin{figure}[h!] 
\begin{center}
\begin{tikzpicture}[scale = 1.8]

\draw[thick] (0,0) rectangle(1,  1);\draw[thick] (1,0) rectangle(2,  1);
\draw[thick] (2,0) rectangle(3,  1); \draw[thick] (3,0) rectangle(4, 1);
\draw[thick] (3,1) rectangle(4,  2); \draw[thick] (4,1) rectangle(5,  2);
\draw[thick] (5,1) rectangle(6,  2); \draw[thick] (6,1) rectangle(7,  2);


\draw (0.5,0.5) node{$1$};
\draw (1.5,0.5) node{$i$};
\draw (2.5,0.5) node{$-1$};
\draw (3.5,0.5) node{$-i$};
\draw (3.5,1.5) node{$-k$};
\draw (4.5,1.5) node{$-j$};
\draw (5.5,1.5) node{$k$};
\draw (6.5,1.5) node{$j$};

\end{tikzpicture}
\caption{ The Eierlegende Wollmilch Sau EW.\\
The gluings of the squares follow the action of the quaternion group.\\
Multiplication by $i$ defines the square on the right, by $j$ the square on top.} \label{fig:EW}
\end{center}
\end{figure}

We construct examples with nice properties arising from exotic translation surfaces. We refer to \cite{fm} for an introduction on translation surfaces and their moduli spaces. 
The  {\it Eierlegende Wollmilch Sau}  and the {\it Ornithorynque} are square-tiled surfaces with remarquable properties. The Eierlegende Wollmilch Sau is described in Figure \ref{fig:EW}.
In both cases, the Veech group of the translation surface is equal to 
$SL(2, \GZ)$.  We consider $\phi$ a pseudo-Anosov homeomorphism in the affine group. When $\phi$ fixes a separatrix, using Veech's zippered rectangles method (see \cite{ve} or \cite{yo}), one can construct a self-similar interval exchange transformation $S$ acting on a transversal of the expanding foliation of $\phi$. We will call such interval exchange transformation a self-similar EW-interval exchange transformation (resp.self-similar Or-interval exchange transformation).

 \begin{theorem}\label{thm:EW}
Let $S$ be a self-similar EW-interval exchange transformation (resp. self-similar Or-interval exchange transformation) and $S'$ an affine interval exchange transformation semi-conjugate to $S$, then $S'$ is topologically conjugate to $S$. In particular, $S'$ has no wandering interval.
 \end{theorem}

 

 {\bf Proof}
 
 We make the proof when the surface is the Eierlegende Wollmilch Sau, denoted by EW. It is exactly similar for the Ornithorynque. 

The interval exchange  $S$ is a 9-interval exchange transformation defined on an interval $I$. The intervals of continuity of $T$ are denoted by 
$I_1, \dots, I_9$.
We consider the function $f$ from $I$ to $\GR$ with value at $x$, the logarithm of the slope of the affine interval exchange transformation $T'$ at $x$. 
The function $f$ is constant on the intervals of continuity of $T$ and is orthogonal to the lengths vector of the intervals of $S$ (see \cite{cage}). In other words, 
$f$ can be also considered as a vector belonging to a co-dimension one subset of $\GR^9$. 
It is enough to prove that, for every $x$ there exist two sequences of positive integers $(n_i)$ and ($m_i)$ and a constant $C$ such that 
\begin{equation}\label{eq:boundedsums}\sum_{j=0}^{n_i}\vert f(S^j(x)\vert < C \textrm{ and } \sum_{j=-m_i}^{0}\vert f(S^j(x)\vert < C.
\end{equation}
This trivially implies that
$$\sum_n \textrm{exp}(\sum_{j=0}^n f(S^j(x))$$ and $$\sum_n \textrm{exp}(\sum_{j=-n}^0 f(S^j(x))$$ 
are divergent series which contradicts the existence of wandering intervals and prove the topological conjugacy between $S$ and $S'$ (see \cite{cage} or \cite{co}). 

We now give some geometric background. The EW is a surface of genus 3 with 4 singularities. 
We recall that the real dimension of the relative homology with real coefficients $H_1(EW, \Sigma, \GR)$ is 9. That's related to the fact that $S$ is a 9-interval exchange transformation. 
The subspace,  $H_1^{(0)}(EW, \Sigma, \GR)$, made of  homology classes with zero holonomy  has dimension 7. 
Since this surface is a covering of the torus, the homology of the surface splits over $\GQ$ into the standard part that comes from the homology of the torus and the space $H_1^{(0)}(EW, \Sigma, \GR)$.
The action of the affine group on $H_1^{(0)}(EW, \Sigma, \GR)$  is an action by a finite group (see \cite{mayo} for precise definitions and a very detailed combinatorial approach to this problem).
let $\phi$ be the pseudo-Anosov homeomorphism from which $S$ is built.The result of Matheus and Yoccoz implies that  $\phi$ acts on $H_1^{(0)}(EW, \Sigma, \GR)$ as an element of finite order. \medskip

In plain terms, the action of $\phi$ on the relative homology is given by a matrix $B$ (a loop in the Rauzy diagram).
From the work of Zorich \cite{zo} and Forni \cite{fo}, we know that the action of $\phi^*$ is responsible for the deviations of ergodic sums at special times of $T$. By construction $f$ is orthogonal to the Perron-Frobenius eigenvector of $B$. Thus $f$ is contained in the space generated by the contracting eigendirection of $B^t$ and its central space. We can forget the projection of $f$ on the contracting direction since it does not change the study of Birkhoff sums (see \cite{co} for instance). 
Thus, up to taking some power of $\phi$ and considering the projection of $f$ on the central space, we can assume that $B^tf = f$.

We consider the sequence $I^{(p)}$ of  intervals on which the induced map $S^{(p)}$ is similar to $S$. Let $C_1^{p)}, \cdots, C_{9}^{p)}$ be the Rokhlin towers with base the continuity intervals of $I^{(p)}$ and $r_1^{(p)}, \dots  r_{9}^{p)}$ their heights. For each $1 \leq k \leq 9$,
$$C_k^{(p)} = \bigcup_{i=0}^{r_k^{(p)}-1}T^i(I_k^{(p)})$$
where $I_1^{(p)}, \dots, I_{9}^{(p)}$ form a partition of $I^{(p)}$.
The previous geometric discussion means that, there exists a constant $K$, such that, for every $p$ and every $k \in \{1, \dots, 9\}$, for every $y$ in the base of the tower $C_k^{(p)}$
\begin{equation}\label{eq:boundedsum-tower}\vert \sum_{i=0}^{r_k^{(p)}} f(S^iy) \vert = \vert \langle (B^t)^pf, e_k \rangle \vert = \vert \langle f, e_p \rangle \vert< K.
\end{equation}
The above formula comes from a classical argument that can be found, in a more general setting, in Zorich's work.
The inequality  \eqref{eq:boundedsum-tower} proves equation \eqref{eq:boundedsums} for the basis of a Rokhlin tower taking $n_p = r_k^{(p)}$. 
\medskip

Since $S^{(p)}$ is minimal (because $S^{(p)}$ is conjugate to $S$ which is a minimal interval exchange transformation), there exists $\kappa >0$ such that, for each point $x \in I_k^{(p)}$, there exists $j<\kappa$, with 
$(S^{(p)})^j(x)$ belongs to $I_k^{(p)}$. In other words, return times to continuity intervals of $S^{(p)}$ are bounded by a universal constant depending only on $S$. 

Let now $x$ be in $s^i(I_k^{(p)})$. The itinerary of $x$ after time $r_k^{(p)}$ is a succession of towers $C_{1(x)}^{(p)}, \dots, C_{j(x)}^{(p)}$. By the previous remark, there exists $0<j(x)<\kappa$ such that $j(x) = k$. Let $$n_p = r_k^{(p)}-i + r_{1(x)}^{(p)} + \cdots r_{j(x)-1}^{(p)} + i$$
It means that we follow the itinerary of $x$ until it comes back to its initial step, the step $i$ of the Rokhlin tower $k$. Thus the Birkhoff sum of $f$ from time $0$ to time $n_p -1$ is the sum along the towers $C_k^{(p)}, C_{1(x)}^{(p)}, \dots, C_{j(x)}^{(p)}$.
Using bound \eqref{eq:boundedsum-tower}, we have:
$$\sum_{j=0}^{n_p-1}\vert f(T^j(x) \vert < K \kappa.$$
This concludes the proof of Theorem \ref{thm:EW}.

\begin{remark}\label{rsau}
The same results hold for every EW-interval exchange transformation (resp. Or-interval exchange transformation) as soon as the slope of the corresponding linear flow is irrational. The proof goes along the same lines but is a little bit more technical.
\end{remark}

\begin{corollary}\label{csau} Let $L$ be a non-recurrent language, satisfying an unflipped order condition,  such that the language $L'$ of Lemma \ref{nrec3} is $L(S)$ for  any  EW-interval exchange transformation or Or-interval exchange transformation $S$ in Theorem \ref{thm:EW} or Remark \ref{rsau}. Then $L$ is not 
 a  natural coding of any affine interval exchange transformation.\end{corollary}
 {\bf Proof}\\
By Lemma \ref{nrec1}, any non-recurrent orbit of $X_L$ coincides with orbits of $X_{L'}$ on an infinite prefix and an infinite suffix. Thus the estimates made in the proof of Theorem \ref{thm:EW}  on the Birkhoff sums of $S$ contradict the criterion in Theorem \ref{paf}.\qed

\end{document}